\documentclass[a4paper, 12pt]{article}

\usepackage[sort&compress]{natbib}
\bibpunct{(}{)}{;}{a}{}{,} 

\usepackage{amsthm, amsmath, amssymb, mathrsfs, multirow, url, subfigure}
\usepackage{graphicx} 
\usepackage{ifthen} 
\usepackage{amsfonts}
\usepackage[dvipsnames]{xcolor}
\usepackage{fullpage}
\usepackage{xr}
\theoremstyle{plain} 
\newtheorem{theorem}{Theorem}
\newtheorem{corollary}{Corollary}

\theoremstyle{definition}

\theoremstyle{remark} 
\newtheorem{ex}{Example}

\newtheorem*{asmp}{Model Assumptions}

\newcommand{\prob}{\mathsf{P}}

\newcommand{\bin}{{\sf Bin}}
\newcommand{\unif}{{\sf Unif}}
\newcommand{\nm}{{\sf N}}

\newcommand{\chisq}{{\sf ChiSq}}
\newcommand{\bet}{{\sf Beta}}

\newcommand{\RR}{\mathbb{R}}

\newcommand{\XX}{\mathbb{X}}

\newcommand{\UU}{\mathbb{U}}

\newcommand{\TT}{\mathbb{T}}

\newcommand{\G}{\mathscr{G}}
\newcommand{\Gbar}{\overline{\mathscr{G}}}
\newcommand{\Gtilde}{\widetilde{\mathscr{G}}}
\newcommand{\K}{\mathscr{K}}

\newcommand{\eps}{\varepsilon}

\newcommand{\iid}{\overset{\text{\tiny iid}}{\,\sim\,}}

\newcommand{\prior}{\mathsf{Q}}

\newcommand{\action}{\mathbb{A}}

\newcommand{\ratio}{\mathcal{R}}

\newcommand{\cred}{\mathscr{C}}

\newcommand{\lPi}{\rotatebox[origin=c]{180}{$\mathsf{\Pi}$}} 
\newcommand{\uPi}{\mathsf{\Pi}}

\newcommand{\uGamma}{\mathsf{\Gamma}}

\title{Decision-making with possibilistic inferential models\footnote{This is an extended version of the paper \citep{imdec.isipta} presented at the 14th {\em International Symposium on Imprecise Probabilities: Theories and Applications}, July 2025, in Bielefeld, Germany.  That conference paper itself is a substantially updated version of an older working paper \citep{imdec}.}}
\author{Ryan Martin \quad Shih-Ni Prim \quad Jonathan P Williams
}
\date{\today}

\begin{document}

\maketitle 

\begin{abstract}   
Inferential models (IMs) are data-dependent, imprecise-probabilistic structures designed to quantify uncertainty about unknowns. As the name suggests, the focus has been on uncertainty quantification for inference and on its reliability properties in that context.  Focusing on a likelihood-based possibilistic IM formulation, the present paper develops a corresponding framework for decision making, and investigates the decision-theoretic implications of the IM’s reliability guarantees. Here we show that the possibilistic IM’s assessment of an action’s quality, defined by a simple Choquet integral, tends not be too optimistic compared to that of an oracle. This ensures that the IM tends not to favor actions that the oracle doesn’t also favor, hence the IM is also reliable for decision making.  We also establish a complementary, large-sample efficiency result that says the IM's reliability isn't achieved by being grossly conservative. In the special case of equivariant statistical models, further connections can be made between the IM’s and Bayesian's recommended actions, from which certain optimality conclusions can be drawn.

\smallskip

\end{abstract}

\section{Introduction}
\label{S:intro}

In some data analysis applications, the goal is to reliably quantify uncertainty about features of the system or population under investigation.  In others, the primary focus is on decision making.  At least intuitively, quality decisions can be made without reliably inferring features of the system or population; but if one could reliably infer the relevant features of the underlying system, then good---perhaps even optimal---decisions ought to be within reach.  The present paper's goal is to investigate this connection between reliable uncertainty quantification and formal decision-making.  

{\em Inferential models}, or IMs for short, first introduced in \citet{imbook, imbasics}, are data-dependent, (imprecise-) probabilistic structures designed to provide reliable uncertainty quantification about unknowns.  The first version of the IM framework relied heavily on random sets and belief functions; in this paper, however, we'll focus on a more recent and arguably simpler formulation that uses possibility theory and, hence, we'll refer to these as {\em possibilistic IMs} \citep[e.g.,][]{martin.basu, martin.partial2, imreview}.  The specific reliability property that (possibilistic) IMs satisfy is called {\em validity}; see \citep{martin.nonadditive, imchar} and Section~\ref{SS:im} below.  Roughly, an IM is valid if its upper (resp.~lower) probabilities assigned to true (resp.~false) hypotheses don't tend to be too small (resp.~large).  Among other things, validity also implies that IM-derived hypothesis tests and confidence regions exactly---not just asymptotically---control frequentist error rates.  It's precisely the validity property that distinguishes IMs from other brands of (imprecise-) probabilistic statistical inference, including Bayesian, fiducial, and confidence distributions \citep{bernardo.smith.book, xie.singh.2012, hannig.review}, Dempster--Shafer belief functions \citep{dempster2008, shafer1976, wasserman1990b, denoeux.li.2018}, possibility measures \citep{balch2012, hose.hanss.2021, dubois2006}, or more general kinds of imprecise probabilities \citep{augustin.etal.bookchapter, walley1991}.  

The IM developments described above are all focused on statistical inference; the formal decision problem remains to be addressed.  To set the scene, a statistical decision problem comes equipped with a set $\action$ of possible actions, a collection $\TT$ of possible states of the world, and a real-valued loss function $\ell_a(\theta)$ that represents the cost associated with taking action $a \in \action$ when the state of the world is $\theta \in \TT$.  If, in addition, a data-dependent probability distribution supported on $\TT$ is available, then it's common to choose an action that minimizes the corresponding expected loss. Generalizations of this basic framework to cases where uncertainty is quantified via an imprecise probability have been developed, and we'll adopt the same generalized decision-theoretic framework here---that is, we'll rank actions according to a suitable upper expected loss, i.e., a Choquet integral; see Section~\ref{SS:decision}.  Our main focus, however, is on the decision-theoretic implications of the possibilistic IM's form and its validity property.  As stated above, if the IM reliably solves the inference problem, then, in some sense, it should also reliably solve the decision problem.  So, our goal is to determine if and how the IM's risk assessment is reliable and what the implications are for decision rules derived from it. 

Following some background in Section~\ref{S:background} on the statistical problem, possibilistic IMs, and existing perspectives on formal decision-making, we define our proposed IM risk assessment: the Choquet integral of the loss function with respect to the IM's output, which is a possibility measure.  Formal decision-making can then proceed by, say, choosing the action to minimize the IM's risk or upper expected loss. Before getting into the properties of the proposed risk assessments and derived decision rules, we pause to develop some intuition about the advantages of this particular formulation.  That is, given the IM's inherent imprecision, the upper expected loss is generally more conservative than a Bayesian posterior expected loss.  But the goal is to be reliable, not just to be conservative, and we argue in Section~\ref{SS:bayes} that the Bayesian risk assessment is overconfident in the following sense: there exists gambles that are acceptable or even desirable to the Bayesian for given data but systematically lose money.  The IM's risk assessment, however, is apparently conservative in just the right way to avoid this kind of overconfidence.  

The paper's first key result shows that the IM's risk assessment tends not to underestimate that of a so-called quasi-oracle who possesses partial information about the true $\Theta$; since the quasi-oracle's information is ``partial,'' the corresponding risk assessment is slightly more conservative than the oracle's.  This is relevant because Bayesians' overconfidence is the result of a tendency to assign risk smaller than that of an oracle, thereby creating an opportunity for systematic loss.  Theorem~\ref{thm:action} shows, however, that validity implies the IMer's and quasi-oracle's assessments tend not to be drastically inconsistent in the sense that existence of an action that the IMer deems to be relatively good but the quasi-oracle deems to be relatively poor is a rare event with respect to the sampling distribution of the data.  It's in this sense, at least, that the valid possibilistic IM is reliable for the decision problem too. Importantly, the IM's decision-making reliability isn't the result of making grossly conservative assessments.  This point is established by considering a large-sample setting and showing that the IM's risk assessment merges asymptotically with the oracle's assessment, meaning that there's no artificial cushion between the IM and oracle assessments that explains the reliability differences.

An important step in these developments is a characterization of the IM's credal set, i.e., the set of precise probabilities dominated by its upper probability.  Existing results \citep[e.g.,][]{destercke.dubois.2014} imply that the IM's credal set consists of confidence distributions: for each $\alpha \in [0,1]$, each probability measure in the IM's credal set assigns probability at least $1-\alpha$ to the $100(1-\alpha)$\% confidence sets.  For group invariant models (Section~\ref{S:group}), this characterization can be strengthened---the ``maximal element'' in the IM's credal set has a familiar form.  Indeed, Fisher's fiducial argument is ideally suited for these models and the resulting fiducial distribution agrees with the Bayesian posterior distribution based on the right invariant Haar prior; the common distribution returned by these two distinct arguments is a confidence distribution. Then \citet{martin.isipta2023} showed that the ``maximal element'' in the IM's credal set is exactly this Bayes/fiducial distribution. A key decision-theoretic consequence of this connection is Theorem~\ref{thm:risk}, which says that, under certain conditions on the loss function, the IM's risk is minimized at the same action where the Bayes/fiducial distribution's risk is minimized.  Combined with the main result in \citet{taraldsen.lindqvist.2013} about decision rules derived from the fiducial distribution, Theorem~\ref{thm:risk} implies the IM's recommended actions are generally high-quality---i.e., not conservative or inefficient---and, in some cases, optimal.  Therefore, while the IM's validity guarantee requires some degree of conservatism, this does not affect the quality of the IM's recommended actions.  

Aside from the conceptual and theoretical results described above, the paper offers a variety of numerical illustrations, including in Section~\ref{S:applications}, with a variety of loss functions, to showcase how the IM approach works, how it compares to Bayes/fiducial, and how flexible it is.  Section~\ref{S:discuss} offers some concluding remarks and a discussion of a few open problems related to the developments in the present paper.

\section{Background}
\label{S:background}

\subsection{Setup, notation, and perspective}

Let $X$ denote the observable data, taking values in a sample space $\XX$.  It'll often be the case that $X$ has components $X_1,\ldots,X_n$, where each $X_i$ itself might split into an independent--dependent variable pair of arbitrary dimension and form.  When the sample size $n$ is relevant to the discussion, we'll write $X^n = (X_1,\ldots,X_n)$.  As is customary, we write upper-case $X$ for the observable data modeled as a random variable, and lower-case $x$ for a generic value/particular realization of $X$.

Next, a statistical model $\{\prob_\theta: \theta \in \TT\}$, consisting of probability distributions supported on (subsets of) $\XX$, is introduced to quantify the variability or aleatory uncertainty in the observable data $X$.  The probability distributions $\prob_\theta$ have a corresponding density/mass function $p_\theta(x)$ and, when $X=x$ is observed, we'll interpret $\theta \mapsto p_\theta(x)$ as the likelihood function.  The ``true value'' of the unknown model parameter will be denoted by $\Theta$ and this entails an assumption that the model is correctly specified, so that the statement ``$X \sim \prob_\Theta$'' is true.  As is typical in the statistics literature, we will assume throughout that prior information about $\Theta$ is fully vacuous.  The goal, then, is to quantify uncertainty about $\Theta$, given the observation $X=x$. 

There are different perspectives on what ``uncertainty quantification'' means.  On the one hand, some might say uncertainty quantification requires a probability distribution for the unknowns, given observed data.  By far the most common such approach is {\em Bayesian inference}, where a prior distribution $\prior$ for the model parameter $\Theta$ is updated to a posterior distribution $\prior_x$ in light of the observed data $X=x$, using Bayes's formula:
\[ \prior_x(d\theta) \propto p_\theta(x) \, \prior(d\theta). \]
There are other less familiar forms of probabilistic uncertainty quantification, including (generalized) fiducial; see the references cited in Section~\ref{S:intro}.  On the other hand, some might be satisfied with mathematically less formal constructs, such as confidence regions or p-values, whose interpretation is drawn solely from the frequentist sampling distribution properties that the corresponding procedures satisfy.  

Our view is that neither of these brands of uncertainty quantification are fully satisfactory.  First, frequentist solutions have no agreed-upon mathematical framework that describes their fixed-data interpretation and calculus; this lack of guidance has surely contributed to the replication crisis in science.  Second, the Bayesian solution depends on the choice of prior distribution which, in real applications, is artificial at least to some degree.  Practical problems lack the information necessary to determine a prior distribution, so at least some elements of $\prior$ must be concocted by the data analyst.  This makes the desired interpretation of the posterior $\prior_x$ as an ``update of prior information'' unwarranted.  In response to this, the so-called {\em pragmatic Bayesians} interpret the posterior as offering ``quick and dirty confidence'' \citep{fraser2011}, but this interpretation is dubious as well.  Indeed, the {\em false confidence theorem} \citep{balch.martin.ferson.2017} says that Bayes solutions, including those based on default-priors, are unreliable in the sense that they tend to assign high posterior probability to hypotheses that are false.  Further details on the false confidence phenomenon, including examples and insights that highlight the scope and cause of the problem, can be found in \citet{martin.nonadditive, martin.belief2024, imreview}. 
It's these shortcomings of existing statistical frameworks for uncertainty quantification that motivate our alternative perspective described next.


\subsection{Possibilistic inferential models}
\label{SS:im}

The first IM developments relied on random sets and belief functions.  More recent developments in \citet{martin.partial2}, building on \citet{plausfn, gim}, define a possibilistic IM by applying a version of the probability-to-possibility transform to the model's relative likelihood.  This latter IM construction is the starting point of the present paper. 

The model and observed data $X=x$ determine the relative likelihood
\[ R(x,\theta) = \frac{p_\theta(x)}{\sup_\vartheta p_\vartheta(x)}, \quad \theta \in \TT. \]
We'll implicitly assume here that the denominator is finite for almost all $x$.  As is typical in the literature, we'll also assume that prior information about $\Theta$ is vacuous. 

The relative likelihood defines a possibility contour, i.e., a non-negative function such that $\sup_\theta R(x,\theta) = 1$ for almost all $x$.  This contour determines a possibility measure that can be used for data-driven uncertainty quantification about $\Theta$, which has been extensively studied in the statistic literature \citep[e.g.,][]{shafer1982, wasserman1990b, denoeux2006, denoeux2014}.  This likelihood-driven possibility has a number of desirable properties.  What it lacks, however, is a justification for why the ``possibilities'' assigned to hypotheses about $\Theta$ are meaningful to determine or at least inform the data analyst's beliefs.  Indeed, with vacuous prior information, there's no Bayesian justification behind these possibility assignments, so justification may only come from a validity-like frequentist reliability guarantee.  But the likelihood-based possibility assignment falls short of this goal.  So, while the relative likelihood provides a useful, data-driven parameter ranking, our conclusion is that $R$ alone is insufficient for reliable statistical inference.  


Fortunately, it's straightforward to achieve the desired goal by ``validifying'' \citep{martin.partial2} the relative likelihood.  This amounts to applying a version of the probability-to-possibility transform \citep[e.g.,][]{dubois.etal.2004, hose2022thesis}, and the result is the possibilistic IM's contour function.  The most basic version is:
\begin{equation}
\label{eq:contour}
\pi_{x}(\theta) = \prob_\theta\bigl\{ R(X,\theta) \leq R(x, \theta) \bigr\}, \quad \theta \in \TT.
\end{equation}
Then the corresponding possibility measure, or upper probability, is defined as 
\begin{equation}
\label{eq:maxitive}
\uPi_{x}(H) = \sup_{\theta \in H} \pi_{x}(\theta), \quad H \subseteq \TT.
\end{equation}
There's a corresponding necessity measure, or lower probability, defined as $\lPi_{x}(H) = 1 - \uPi_{x}(H^c)$, but this plays a very limited role here in this paper and in our approach to statistical inference more generally.  Some brief additional points along these lines are given at the end of the present subsection. 

As in Section~\ref{S:intro}, an essential feature of this IM construction is its {\em validity} property:
\begin{equation}
\label{eq:valid}
\sup_{\theta \in \TT} \prob_\theta\bigl\{ \pi_{X}(\theta) \leq \alpha \bigr\} \leq \alpha, \quad \text{$\alpha \in [0,1]$}. 
\end{equation}
Property \eqref{eq:valid} has a number of important consequences.  First, it immediately implies that 
\begin{equation}
\label{eq:region}
C_\alpha(x) = \{\theta \in \TT: \pi_{x}(\theta) > \alpha\}, \quad \alpha \in [0,1]
\end{equation}
is a $100(1-\alpha)$\% frequentist confidence set in the sense that 
\[ \sup_{\theta \in \TT} \prob_\theta\bigl\{ C_\alpha(X) \not\ni \theta \bigr\} \leq \alpha, \quad \alpha \in [0,1]. \]
Second, from \eqref{eq:maxitive} and \eqref{eq:valid}, 
\begin{equation}
\label{eq:valid.alt}
\sup_{\theta \in H} \prob_\theta\bigl\{ \uPi_{X}(H) \leq \alpha \bigr\} \leq \alpha, \quad \text{$\alpha \in [0,1]$, $H \subseteq \TT$}. 
\end{equation}
In words, a valid IM assigns possibility $\leq \alpha$ to true hypotheses at rate $\leq \alpha$ as a function of data $X$.  This gives the IM its ``inferential weight''---\eqref{eq:valid.alt} implies that $\uPi_{x}(H)$ is not expected to be small when $H$ is true, so one is inclined to doubt the truthfulness of a hypothesis $H$ if $\uPi_{x}(H)$ is small.  Third, the above property ensures that the possibilistic IM is safe from false confidence \citep{balch.martin.ferson.2017, martin.nonadditive}, unlike default-prior Bayes and fiducial solutions.  
Further properties of possibilistic IMs are discussed in \citet{imreview} and in the sections that follow. 

The IM output is a coherent imprecise probability and, consequently, there's a non-empty set of precise probabilities that are compatible with it.  This is called the {\em credal set} associated with the possibilistic IM's output $\uPi_x$, and is expressed mathematically as 
\begin{equation}
\label{eq:credal}
\cred(\uPi_x) = \{\prior_x \in \text{probs}(\TT): \prior_x(\cdot) \leq \uPi_x(\cdot)\}, 
\end{equation}
where $\text{probs}(\TT)$ is the set of countably additive probabilities supported on (the Borel $\sigma$-algebra of subsets of) $\TT$ and ``$\prior_x(\cdot) \leq \uPi_x(\cdot)$'' means the inequality holds for all measurable events on $\TT$.  Of course, the members of $\cred(\uPi_x)$ depend on $x$ because $\uPi_x$ does, but these may not correspond to Bayesian posterior distributions under any prior.  Fortunately, an interpretation can be given to the members of $\cred(\uPi_x)$ thanks to a well-known characterization \citep[e.g.,][]{destercke.dubois.2014, cuoso.etal.2001}: 
\[ \prior_x \in \cred(\uPi_x) \iff \text{$\prior_x\{C_\alpha(x)\} \geq 1-\alpha$, all $\alpha \in [0,1]$}, \]
where $C_\alpha(x)$ is as defined in \eqref{eq:region} with $\pi_x$ the contour corresponding to $\uPi_x$.  Since $C_\alpha(x)$ is a $100(1-\alpha)$\% confidence set, and the elements $\prior_x$ in the credal set assign at least probability $1-\alpha$ to $C_\alpha(x)$, there's good reason to call these elements {\em confidence distributions}.  This definition of confidence distributions agrees with that given in \citet{taraldsen.2021.cd} and generalizes those commonly found in the literature \citep[e.g.,][]{xie.singh.2012, schweder.hjort.cd.discuss, prsa.conf}.  This credal set characterization and its interpretation as a collection of confidence distributions is important in what follows. 

We end this subsection with two technical remarks.  The first is about why the IM's lower probability $\lPi_x$ plays effectively no role in our analysis.  As is customary, we adopt a Popperian falsificationist perspective \citep[e.g.,][]{popper1959, popper1962} wherein scientific theories can be refuted based on an incompatibility with data but can never be confirmed.  For instance, p-values correspond to valid possibilistic IM upper probabilities \citep[e.g.,][App.~G]{imreview}, and one can reject a null hypothesis $H_0$ when $\uPi_x(H_0)$ is small but cannot, under any circumstances, accept $H_0$.  If there is no threshold $\beta \in (0,1)$ such that $\lPi_x(H_0)$ exceeding $\beta$ would warrant an ``accept $H_0$'' conclusion, then it must be that $\lPi_x$ carries no weight on its own.  So, since statistical inference is about refutation, and refutation is based on demonstrating incompatible or implausibility, it is without loss of flexibility or generality that we focus exclusively on the IM's upper probability output $\uPi_x$.  

Second, while the focus here is on inference and decision-making in the context of inference on a model parameter about which prior information is vacuous, there are other contexts in which the IM formulation can be applied.  Prediction of future observations is one such context, and there are some nontrivial differences in the technical details.  For instance, there are actually two distinct notions of the validity property described above, called {\em weak} and {\em strong} validity.  These two are equivalent in the context of the present paper, hence we only need to define one such notion.  But it turns out, for instance, that the now widely-studied conformal prediction framework \citep[e.g.,][]{vovk.shafer.book1, angel.bates.2023} is a special case of a strongly valid possibilistic IM framework described here; see \citet[][Sec.~6.3]{imreview}.  Further details can be found in \citet{imconformal, imconformal.supervised}, \citet{caprio.etal.isipta25}, \citet{caprio.joys}, and \citet{chau.etal.2026.conformal}.

\subsection{Decision theory}
\label{SS:decision}

Decision theory is a mathematical framework for analyzing the choices agents make in the presence of uncertainty.  One key aspect of this theory is characterizing the ``optimal'' behavior for an agent faced with a decision problem.  This requires comparing the agent's strategies and ranking them by preference.  Early efforts \citep[e.g.,][]{bernoulli.risk} went directly to considerations of loss (or negative utility), ranking strategies according to their expected loss and, hence, defining the optimal strategy as one that minimizes expected loss.  Other authors took a different route, by considering general preference orders on strategies, but reached effectively the same destination.  For example, the celebrated theorem of \citet{vnm.games} says that if the preference order satisfies certain rationality axioms, then there exists a loss function such that the ranking of strategies according to preferences is equivalent to the ranking by expected loss.  Similar conclusions have been reached by \citet{savage1972} and others, hence, the {\em minimize-expected-loss} principle or, equivalently, the {\em maximize-expected-utility} principle.   

For the statistical decision problem, let $(\theta, a) \mapsto \ell_a(\theta)$ denote a loss function that measures the loss incurred by taking an action $a \in \action$ when the ``state of nature'' is $\theta$.  Common examples include squared-error loss, with $\ell_a(\theta) = \|a-\theta\|^2$, and 0--1 loss, with 
\[ \ell_a(\theta) = 1(\theta \in H, a=1) + 1(\theta \not\in H, a=0), \]
where $H \subset \Theta$ is a (null) hypothesis, and $a=1$ and $a=0$ denote ``reject'' and ``do not reject'' $H$, respectively.  We will assume throughout that the loss is non-negative.  While negative losses might make sense in some contexts, these gains would typically be bounded and so the loss could be made non-negative by adding a constant.  This constant shift won't affect judgments of the relative quality of actions. 

Of course, an oracle who knows the true $\Theta$ could easily compare possible actions according to their corresponding loss values, $\ell_a(\Theta)$.  We'll refer to $a \mapsto \ell_a(\Theta)$ as the {\em oracle's assessment} of action $a$.  Naturally, the oracle would take the action $a^\star = \arg\min_a \ell_a(\Theta)$ and incur minimal loss.  But none of us have oracle powers, so $\ell_a(\Theta)$ is out of our reach and, therefore, a different approach is required.  When uncertainty about $\Theta$, given $X=x$, is quantified via a probability $\prior_x$, like a Bayesian posterior distribution, a typical strategy is to rank the candidate actions according to their expected loss, $a \mapsto \prior_x \ell_a$, where $\prior_x f = \int f(\theta) \, \prior_x(d\theta)$ is the ordinary $\prior_x$-expected value of $f$.  Then, mimicking the oracle, one minimizes the expected loss, $\hat a(x) = \arg\min_a \prior_x \ell_a$. 

Other authors have argued that requiring a single precise probability to quantify uncertainty puts too much of a burden on the data analyst, and that a decision-theoretic framework based on more general imprecise probabilities is better suited for practical applications.  Excellent reviews of decision theory from an imprecise probability perspective can be found in \citet{huntley.etal.decision} and \citet{denoeux.decision.2019}.  Roughly, given a loss function like above, if uncertainty is quantified via an imprecise probability, then it's only natural to extend the minimize-expected-loss principle by replacing the expected loss, $\prior_x \ell_a$, with an appropriate generalization, and then optimizing to find the best action \citep[e.g.,][]{gilboa1987, gilboa.schmeidler.1989}.  Indeed, it was shown in \citet{gilboa.schmeidler.1994} that an appropriate generalization of the expected loss is obtained via the Choquet integral; see Appendix~C of \citet{lower.previsions.book}, Section~4.1 in \citet{denoeux.decision.2019}, and Section~\ref{S:decisions} below for the special case of possibility measures. 

Regardless of the particular imprecise probabilistic framework one is working in, there are two versions of ``expected loss'' that can be considered: a lower and an upper expected loss.  From here, there are numerous criteria that generalized the basic minimize-expected-loss principle.  One of the most common of these generalizations is the criterion that chooses the action that minimizes the upper expected loss.  For obvious reasons, this strategy is commonly referred to as {\em minimax}, and there is extensive work about this in the literature on imprecise probability \citep{huntley.etal.decision} and in robust Bayesian inference \citep{berger1984, berger1985, vidakovic2000}. 
Von Neumann and Morgenstern-style axiomatic justification for this minimax criterion can be found in \citet{jaffray1988, jaffray1989}, \citet{gilboa1987}, and \citet{schmeidler.1989}.  A decision-theoretic framework based on upper probabilities, derived from the {\em model-free generalized fiducial} argument developed in \citet{williams.cpfid} is presented in \citet{williams2024decision}. Moreover, special decision-theoretic attention has been paid to ``partially consonant'' belief models \citep[e.g.,][]{walley1987, giang.shenoy.2011}, which include the possibility measures considered here.  Beyond the somewhat pessimistic minimax criterion, there's the more optimistic {\em maximin} criterion that maximizes the lower expected loss and the even more optimistic {\em minimin} criterion that minimizes the lower expected loss.  The so-called {\em Hurwicz criterion} balances the pessimism and the optimism by choosing the action that minimizes a convex combination of the upper and lower expected loss.  See \citet{denoeux.decision.2019} for a detailed discussion of these and other decision criteria.  For reasons given in Section~\ref{SS:minimax}, we focus on the minimax strategy.  


\section{IMs and decision-making}
\label{S:decisions}

\subsection{IM-based risk assessments}
\label{SS:imrisk}

We start by developing the decision-theoretic framework for possibilistic IMs.  When the upper probability $\uPi_x$ is interpreted as an upper envelope on a collection of ordinary probabilities, then, naturally, its extension to an upper expectation/prevision is  
\begin{equation}
\label{eq:upper.expectation}
\uPi_x f = \sup\{ \prior_x f: \prior_x \in \cred(\uPi_x)\}, 
\end{equation}
where $\cred(\uPi_x)$ is the credal set \eqref{eq:credal} and $f$ is a suitable real-valued function defined on $\TT$.  In light of the discussion in Section~\ref{SS:im}, the upper envelope has some further intuition, that is, $\uPi_x f$ is largest of the ordinary expected values, $\prior_x f$, corresponding to confidence distributions $\prior_x$ compatible with $\uPi_x$.  There is an associated lower expectation, $\lPi_x f = -\uPi_x(-f)$, but, for reasons discussed in Section~\ref{SS:im} and \ref{SS:minimax}, this won't be needed. 

The above extension of an upper probability to an upper expectation can be equivalently described as a {\em Choquet integral} \citep[e.g.,][]{choquet1953, huber1973.capacity, denneberg1994}.  The focus here is on IMs whose upper probability $\uPi_x$ is a possibility measure.  For these, Propositions~7.14 and 15.42 in \citep{lower.previsions.book}, in increasing generality, establish that the Choquet integral of a non-negative function $f: \TT \to [0,\infty)$ with respect to the possibility measure $\uPi_x$, if it exists, is given by 
\begin{equation}
\label{eq:choquet}
    \uPi_x f = \int_0^1 \Bigl\{ \sup_{\theta: \pi_x(\theta) > s} f(\theta) \Bigr\}  \, ds, 
\end{equation}
where $\pi_x$ is the possibility contour corresponding to $\uPi_x$.  (The integral in the above display is a Riemann integral, which exists because the integrand is monotone in $s$.) Existence of the Choquet integral requires that $f$ be ``previsible'' with respect to $\uPi_x$ \citep[][Def.~15.6]{lower.previsions.book}, a condition we'll silently assume about the loss functions below.
Importantly, the Choquet integral determines the upper expectation, hence the right-hand sides of \eqref{eq:upper.expectation} and \eqref{eq:choquet} are equal---this implies that evaluation of $\uPi_x f$ doesn't require optimization over distributions $\prior_x$ in the credal set $\cred(\uPi_x)$!  

Since the loss function is assumed to be non-negative, either of the above two equivalent formulas can be applied directly to define an upper risk/expected loss,  
\begin{equation}
\label{eq:upper.risk}
a \mapsto \uPi_x \ell_a. 
\end{equation}
Assuming $\uPi_x \ell_a$ is finite at least for some actions $a$, this can be used to assess the quality of different actions, relative to the given loss function and the IM's data-dependent possibility measure.  In particular, we propose to adopt the aforementioned minimax criterion, which defines the IMer's $x$-dependent ``optimal'' action to be 
\begin{equation}
\label{eq:best}
\hat a(x) = \arg\min_a \uPi_x \ell_a. 
\end{equation}
Section~\ref{SS:examples1} below presents a few illustrative examples of this proposal and Section~\ref{SS:minimax} offers our justification for focusing on the minimax criterion.


\subsection{Examples}
\label{SS:examples1}

\begin{ex}
\label{ex:test}
Consider a hypothesis testing problem, with $H_0: \Theta \in \TT_0$ and $H_1: \Theta \in \TT_1$, where $\TT_0 \cap \TT_1 = \varnothing$.  The action space is $\action = \{0,1\}$ and the standard 0--1 loss is 
\[ \ell_a(\theta) = 1(\theta \in \TT_0, \, a=1) + 1(\theta \in \TT_1, \, a=0), \quad \theta \in \TT_0 \cup \TT_1. \]
Then it's easy to see that the possibilistic IM's assessment is 
\[ \uPi_x \ell_a = \uPi_x(\TT_1) \, 1(a=0) + \uPi_x(\TT_0) \, 1(a=1), \]
and the corresponding IM's minimax optimal action is 
\[ \hat a(x) = \begin{cases} 0 & \text{if $\uPi_x(\TT_0) < \uPi_x(\TT_1)$} \\ 1 & \text{if $\uPi_x(\TT_1) < \uPi_x(\TT_0)$}. \end{cases} \]
In the event that $\uPi_x(\TT_0)=\uPi_x(\TT_1)$, one might set $\hat a(x) = \{0,1\}$ or flip a coin to decide.  
\end{ex}

\begin{ex}
\label{ex:location}
Consider a simple location model, on $\XX=\RR$, where the density $p_\theta$ corresponding to the distribution $\prob_\theta$ is of the form $p_\theta(x) = p(x - \theta)$, for some density $p$ symmetric around 0, e.g., Gaussian, Student-t, etc.  We'll assume that $n=1$, but neither this nor symmetry of $p$ are absolutely necessary; this just makes the calculations easy to do by hand and helps reveal some interesting structure investigated in Section~\ref{S:group}.  

In this case, the relative likelihood is $R(x,\theta) = p(x-\theta) / p(0) \propto p(x-\theta)$.  Also, since $p$ is symmetric about 0, the IM contour for $\Theta$, as defined in \eqref{eq:contour}, reduces to 
\begin{align*}
\pi_x(\theta) & = \prob_\theta\{ p(X-\theta) \leq p(x-\theta) \} \\
& = \prob_\theta\{ |X-\theta| \geq |x-\theta| \} \\
& = 2\{1 - P(|x-\theta|)\}, \qquad \theta \in \RR,
\end{align*}
where $P$ is the cumulative distribution function corresponding to the density $p$.  Then the $100(1-\alpha)$\% confidence interval in \eqref{eq:region} is 
\[ C_\alpha(x) = \bigl[ x - P^{-1}(1-\tfrac{\alpha}{2}), \, x + P^{-1}(1-\tfrac{\alpha}{2}) \bigr]. \]
To define a decision problem, consider ``point estimation'' where the action space $\action$ is the parameter space.  As is customary, we'll use the squared error loss, i.e., $\ell_a(\theta) = (a - \theta)^2$.  It'll be clear in what follows that $\uPi_x \ell_a$ is well-defined only if the density $p$ has sufficiently thin tails that its variance $v(p)$ is finite, so we assume this throughout.  

\begin{figure}[t]
\begin{center}
    \scalebox{0.65}{\includegraphics{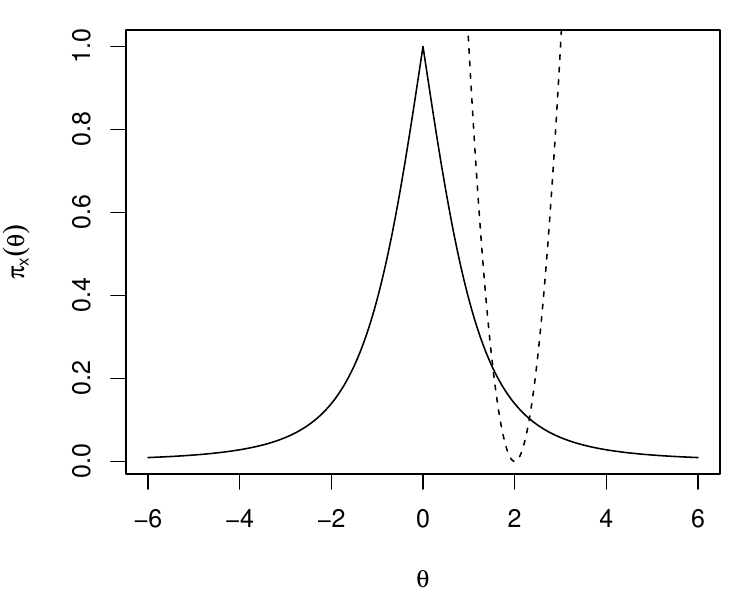}}
\end{center}
\caption{Plot of $\pi_x(\theta)$, for $x=0$, when $p$ is a Student-t density with 3 degrees of freedom; dashed line shows $\ell_a(\theta)$, when $a=2$.}
\label{fig:pl.loss.1}
\end{figure}

Figure~\ref{fig:pl.loss.1} shows both the contour function $\pi_x(\theta)$ and the loss function $\ell_a(\theta)$ for specific values of $x$ and $a$, with $p$ a Student-t density.  Observe that, on any level set of $\pi_x$, the loss function attains its maximum at the boundary.  With this observation, the Choquet integral formula in \eqref{eq:choquet} gives 
\begin{align*}
\uPi_x \ell_a & = \int_0^1 \sup_{\theta: 2\{1 - P(|x-\theta|)\} > s} (a - \theta)^2 \, ds \\
& = \int_0^1 \max_{+,-} \bigl\{x - a \pm P^{-1}(1-\tfrac{s}{2}) \bigr\}^2 \, ds \\
& = \int_0^1 \bigl\{ |x-a| + P^{-1}(1-\tfrac{s}{2}) \bigr\}^2 \, ds \\
& = \int_{-\infty}^\infty \bigl( |x-a| + |z| \bigr)^2 \, p(z) \, dz \\
& = |x-a|^2 + 2 \, m(p) \, |x-a| + v(p), 
\end{align*}
where $m(p) = \int |z| \, p(z) \, dz$ and $v(p) = \int z^2 \, p(z) \, dz$. This is clearly minimized at $\hat a(x) = x$, which is the Bayes rule and Pitman's optimal equivariant estimator \citep{pitman1939a}.  This connection between the IM risk minimizer and the Bayes rule in the location parameter problem is not a coincidence, as we show in Section~\ref{S:group}.  Importantly, while the Bayes and IM risk assessments have the same minimizer, they are not the same function.  Indeed, $\uPi_x \ell_a$ is a quadratic function in $|x-a|$ whereas the Bayes risk $\prior_x \ell_a$ is quadratic in $a$. This implies that the IM risk exceeds the Bayes risk, which is a general phenomenon and a key feature, as discussed in Section~\ref{SS:bayes}.   
\end{ex}

\begin{figure}
    \centering
    \subfigure[Contour function]{\includegraphics[width=0.49\linewidth]{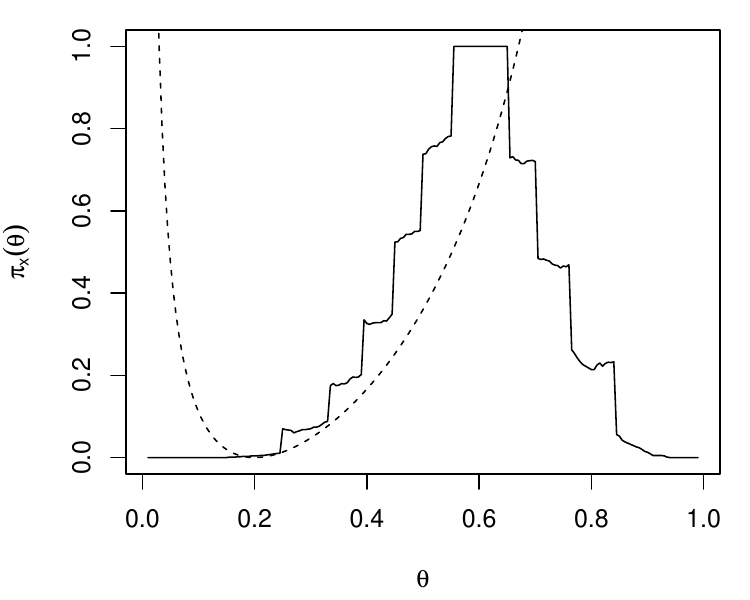}}
    \subfigure[Risk assessment]{\includegraphics[width=0.49\linewidth]{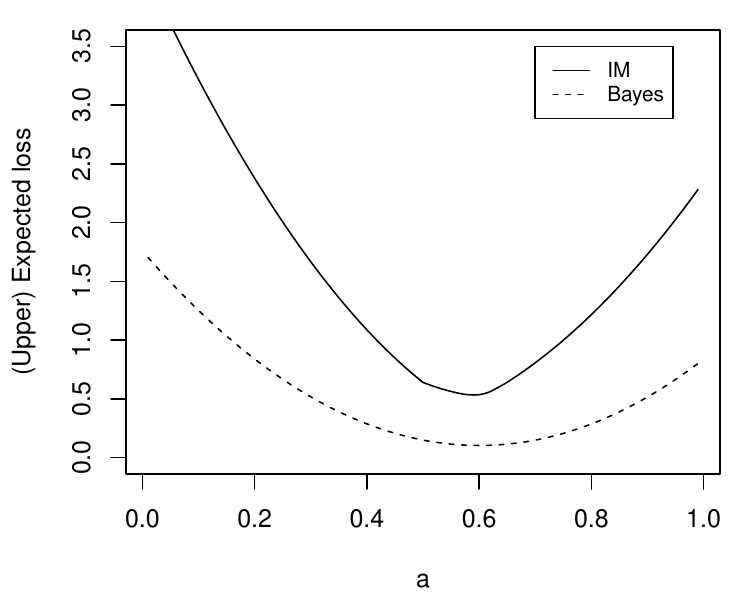}}
    \caption{Results from Example~\ref{ex:binomial}. Panel~(a) shows the IM contour function $\pi_x(\theta)$ (solid) based on $n=10$ and $x=6$; dashed line shows $\theta \mapsto \ell_a(\theta)$, when $a=0.2$. Panel~(b) shows the Bayes (dashed) and IM (solid) risk assessments.}
    \label{fig:binomial}
\end{figure}

\begin{ex} 
\label{ex:binomial}
Let $X \sim \bin(n, \theta)$.  The standard, default-prior Bayesian solution---which appears in the original essay of \citet{bayes1763}---assumes a uniform prior, i.e., 
\[ \Theta \sim \unif(0,1) \equiv \bet(1, 1). \]
This is a conjugate prior, so the posterior distribution of $\Theta$, given $X=x$, is 
\[ (\Theta \mid x) \sim \prior_x := \bet(x + 1, n - x + 1). \]
For the possibilistic IM solution, the relative likelihood is 
\[
R(x, \theta) = \Bigl( \frac{n\theta}{x} \Bigr)^{x} \Bigl( \frac{n(1-\theta)}{n-x} \Bigr)^{n-x}.
\]
The corresponding possibility contour $\pi_x(\theta)$ as defined in \eqref{eq:contour} can be exactly evaluated numerically---no Monte Carlo needed---but there's no closed-form expression.  Likewise, the Choquet integral must be calculated numerically.  On that note, now let's consider a weighted version of the squared error loss: 
\[
    \ell_a(\theta) = \frac{(a-\theta)^2}{\theta(1-\theta)}, \theta \in [0, 1].
\]
We carried out an experiment with $n=10$. Figure~\ref{fig:binomial}(a) shows the contour function with $\Theta=0.5$ and the loss function with one specific choice of action. The contour function shows some steps, due to the discrete nature of the binomial distribution. Figure~\ref{fig:binomial}(b) shows the Bayes and IM risk assessments as functions of the action $a$. Note that the two curves have approximately the same minimizers, and yet, again, the Bayes assessment is uniformly smaller than the IM assessment; see Section~\ref{SS:bayes}.
\end{ex}

\subsection{Why minimax?}
\label{SS:minimax}

Our proposed decision-making framework is based on an adoption of the minimax criterion.  After showing a few examples of the minimax strategy in action, we're ready to offer a justification for that choice.  

If one accepts our philosophical explanation in Section~\ref{SS:im} of why the IM's lower probability $\lPi_x$ has little/nothing to contribute to the inference problem, then one must also accept that the corresponding lower expectation $\lPi_x \ell_a$ has little/nothing to contribute to the decision problem.  If sole focus lands on the upper expected loss $a \mapsto \uPi_x\ell_a$, then choosing $\hat a$ as the corresponding minimizer is justified not only by the axiomatic developments surveyed in \citet{denoeux.decision.2019} but by common sense.  Indeed, we have no ordering among the members $\prior_x$ of the credal set $\cred(\uPi_x)$ that says one is ``more likely'' than another, so the only way to ensure that $\prior_x \ell_a$ is small is to make the supremum in \eqref{eq:upper.expectation} small, which amounts to minimizing the upper bound $\uPi_x \ell_a$.  

To better understand our preference for minimizing the upper expected loss, it may help to consider some alternatives.  These, of course, will involve consideration of the lower expected loss $\lPi_x \ell_a$.  Since the loss is non-negative, it's easy to check that 
\[ \lPi_x\ell_a = \int_0^1 \Bigl\{ \inf_{\theta: \pi_x(\theta) > s} \ell_a(\theta) \Bigr\} \, ds. \]
Then, say, the minimin criterion corresponds to minimizing $a \mapsto \lPi_x \ell_a$, and it's likewise easy to see that this minimizer boils down to the {\em plug-in action} 
\[ \hat a^\text{\sc pi}(x) = \arg\min_a \lPi_x \ell_a = \arg\min_a \ell_a(\hat\theta_x), \]
where $\hat\theta_x$ is the maximum likelihood estimator.  The high degree of optimism is apparent in this expression---it just mimics what an oracle who knows $\Theta$ would do, pretending that $\hat\theta_x$ is the true $\Theta$.  Note that while the risk assessment associated with the plug-in action is overly optimistic, the action itself may be quite reasonable.  For instance, in Examples~\ref{ex:test}--\ref{ex:location} in Section~\ref{SS:examples1}, the actions based on the minimax and this new minimin criteria are actually the same; therefore, while the risk assessments might be very different, there's no sense in which the minimin action is ``too optimistic'' and the minimax action is ``too pessimistic.''  When the two suggested actions are different, however, we think it's fair to characterize the minimin action as ``naive,'' since it fails to make full use of the uncertainty quantification provided by the IM output.  More generally, one might adopt a Hurwicz criterion, which amounts to minimizing 
\[ a \mapsto \varrho \, \uPi_x \ell_a + (1-\varrho) \, \lPi_x \ell_a, \]
where $\varrho \in [0,1]$ is a pessimism index.  As explained above, we have no justification for being less than fully pessimistic, which suggests taking $\varrho=1$, yielding the minimax action again.  And if we take $\varrho < 1$, then either (a)~we get back the same minimax action because the minimax and minimin actions agree, or (b)~we shrink the minimax action $\hat a$ towards the naive, minimin plug-in action. It's not clear why shrinking toward the minimin action would be desirable if the minimin action itself is ``naive.''

\subsection{Comparison with Bayes assessments}
\label{SS:bayes}


The two previous examples in Section \ref{SS:examples1} reveal that there's a non-negligible difference between the Bayes and IM risk assessments, and the goal of the present section is to flesh out these differences more thoroughly and precisely---and to explain the implications of these differences.  In particular, the point that we wish to make here is twofold.  First, we argue that the {\em Bayesian assessments are overconfident}, thereby putting the decision-maker at risk; while we stop short here of stating and proving a formal theorem, we're shining a light on a decision-making version of the false confidence phenomenon.  The second point is that the IM assessments are not overconfident in the sense that Bayesian assessments are overconfident.  To be clear, this argument isn't specific to Bayesian decision-making, but probabilistic decision-making more generally; that is, the points below also apply to assessments that involve integrating the loss function with respect to any brand of fiducial or confidence distribution for $\Theta$.    

The jumping off point is the realization that Bayesian/probabilistic risk assessments, $a \mapsto \prior_x \ell_a$, are designed to mimic the oracle assessment, $a \mapsto \ell_a(\Theta)$.
That is, in the unrealistically ideal case when the posterior $\prior_x$ precisely knows the true $\Theta$, as often happens in the large-sample limit, the Bayesian assessment equals the oracle's.  In every realistic case, however, the expected loss with respect to the posterior distribution will return values bigger or smaller than $\ell_a(\Theta)$, depending on $a$.  Cases where $\prior_x \ell_a > \ell_a(\Theta)$ aren't concerning, but cases where $\prior_x \ell_a < \ell_a(\Theta)$ are problematic---even dangerous---in the following sense.  According to de Finetti, the posterior expectation can be interpreted as the Bayesian's ``fair price'' for the gamble in question, given data.  In our present context, since $\ell_a(\theta)$ is a loss, i.e., negative payoff, $\prior_x \ell_a$ is most naturally interpreted as the greatest lower bound on how much he's willing to accept as payment for taking action $a$, given $X=x$.  That is, the Bayesian is willing to accept $\prior_x \ell_a + \eps$ dollars to take action $a$ since his net payoff is then $\prior_x \ell_a + \eps - \ell_a(\Theta)$ dollars and, based on his quantification of uncertainty about $\Theta$, given $X=x$, he expects this net payoff to be positive.  Therefore, under these circumstances, the Bayesian views taking action $a$ in exchange for a payment that exceeds $\prior_x \ell_a$ as a desirable transaction.  We will demonstrate below that, despite these being desirable actions to the Bayesian, there exists actions such that the above net payoff tends to be negative.  This is problematic because it means a (clever or just lucky) adversary can offer up transactions that are judged to be desirable to our Bayesian but will put him at risk of systematically losing money---at least more money than the inherently more conservative IM assessment.  A couple illustrations follow. 

First, a bit more formality.  For a data-dependent action $a(x)$, define the Bayesian's $\eps$-compensated net payoff/winnings as 
\[ W_\eps^\text{\sc ba}(x; \Theta) = \prior_x \ell_{a(x)} + \eps - \ell_{a(x)}(\Theta). \]
As explained above, this ``gamble'' (treated as a function of $\Theta$ for fixed $x$) is {\em desirable} to the Bayesian because his posterior expected payoff with respect to $\prior_x$ is positive.  Similarly, the IMer's $\eps$-compensated net payoff/winnings 
\[ W_\eps^\text{\sc im}(x; \Theta) = \uPi_x \ell_{a(x)} + \eps - \ell_{a(x)}(\Theta), \]
and this gamble is {\em not-undesirable} to the IMer: her upper expected payoff with respect to $\uPi_x$ is positive.  For simplicity, we'll take $\eps = 0$ and say that the corresponding gambles above are {\em acceptable} and {\em not-unacceptable} to the Bayesian and IMer, respectively.  

\begin{figure}[t]
\begin{center}
\scalebox{0.65}{\includegraphics{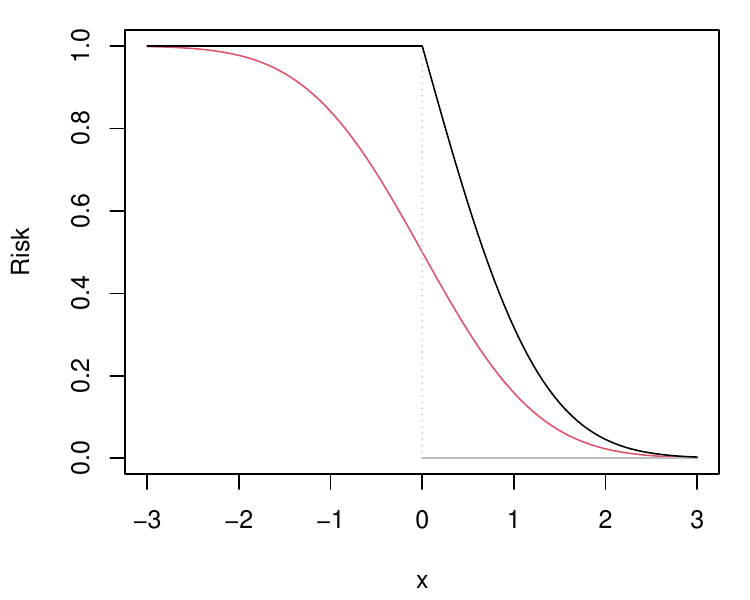}}
\end{center}
\caption{Plots of the Bayes (red) and IM (black) risks, $\prior_x(\TT_0)$ and $\uPi_x(\TT_0)$, associated with the Gaussian test of $H_0: \Theta \in \TT_0 = (-\infty, 0]$ and action $a=1$ (``reject $H_0$''), as a function of data $x$.  The oracle risk (gray) is the indicator function of $\TT_0$.}
\label{fig:gauss.test}
\end{figure}

As a first concrete illustration, consider the hypothesis testing scenario in Example~\ref{ex:test} above, with zero--one loss and the global alternative $\TT_1 = \TT_0^c$.  Consider the action $a=1$, which corresponds to ``reject the null.''  The Bayesian would be willing to take action $a=1$ in exchange for at least $\prior_x(\TT_0)$ dollars, so the following gamble is acceptable:
\[ W_0^\text{\sc ba}(X; \Theta) = \prior_X(\TT_0) - 1(\Theta \in \TT_0). \]
The fact that this gamble is acceptable/fair implies that $W_0^\text{\sc ba}(X; \Theta) \geq 0$ when $\Theta \in \TT_1$ and $W_0^\text{\sc ba}(X; \Theta) \leq 0$ when $\Theta \in \TT_0$; this is effectively the Bayesian's no-sure-loss property.  But note that, except in the most extreme circumstances, the realized payoff when $\Theta \in \TT_0$ and action $a=1$ is wrong is {\em strictly negative}---this is because the posterior probability assigned to the correct $\TT_0$ will surely be strictly less than 1.  In contrast, the IMer would be willing to take action $a=1$ in exchange for at least $\uPi_x(\TT_0)$ dollars, so the following gamble is not-unacceptable:
\[  W_0^\text{\sc im}(X; \Theta) = \uPi_X(\TT_0) - 1(\Theta \in \TT_0). \]
Just like with the Bayesian assessment, $W_0^\text{\sc im}(X; \Theta) \geq 0$ when $\Theta \in \TT_1$ and $W_0^\text{\sc im}(X; \Theta) \leq 0$ when $\Theta \in \TT_0$.  A key difference, a consequence of the IM's possibilistic form, is that $\uPi_X(\TT_0) = 1$ exactly whenever the model maximum likelihood estimator, $\hat\theta_X$, is contained in $\TT_0$.  That is, while the Bayesian's net payoff for taking action $a=1$ would generally be negative when $\Theta \in \TT_0$, the IMer's net payoff would generally be non-negative and would often be exactly 0.  For a simple illustration, consider the Gaussian model $X \sim \nm(\Theta,1)$ with $\TT_0 = (-\infty,0]$.  The IM solution was derived above in Example~\ref{ex:location} and the simple, flat-prior Bayes solution has posterior $\prior_x = \nm(x,1)$.  Figure~\ref{fig:gauss.test} plots the Bayes and IM risk assessments associated with action $a=1$, i.e., $\prior_x(\TT_0)$ and $\uPi_x(\TT_0)$, as functions of the data point $x$.  The key observation is that, for a range of $x$ values, the Bayesian's assessment is below the oracle's assessment---a sign of overconfidence, susceptibility to losing money---while the IMer's assessment is always on or above the oracle's.  This is obviously a very simple illustration, but a similar phenomenon occurs in more modern/practical problems, such as logistic regression; see Section~\ref{SS:logistic}.   

\begin{figure}
    \centering
    \subfigure[Risk functions]{\includegraphics[width=0.49\linewidth]{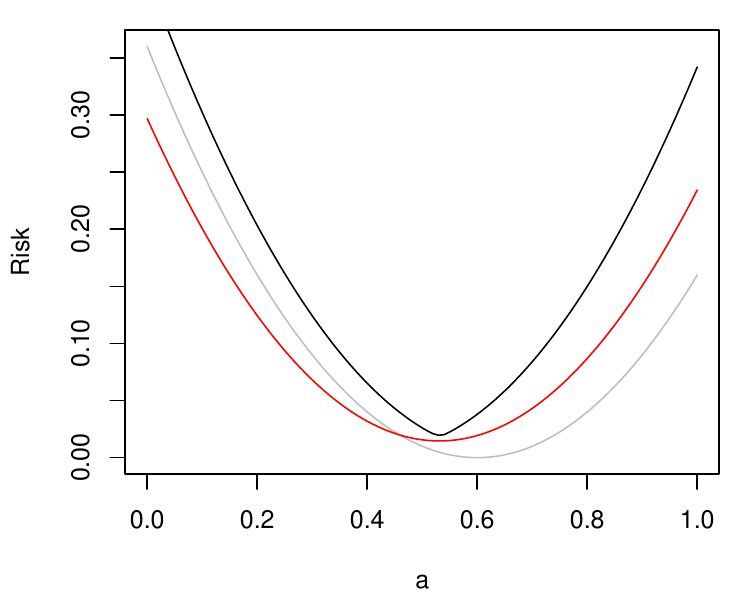}}
    \subfigure[Distribution of net payoff]{\includegraphics[width=0.49\linewidth]{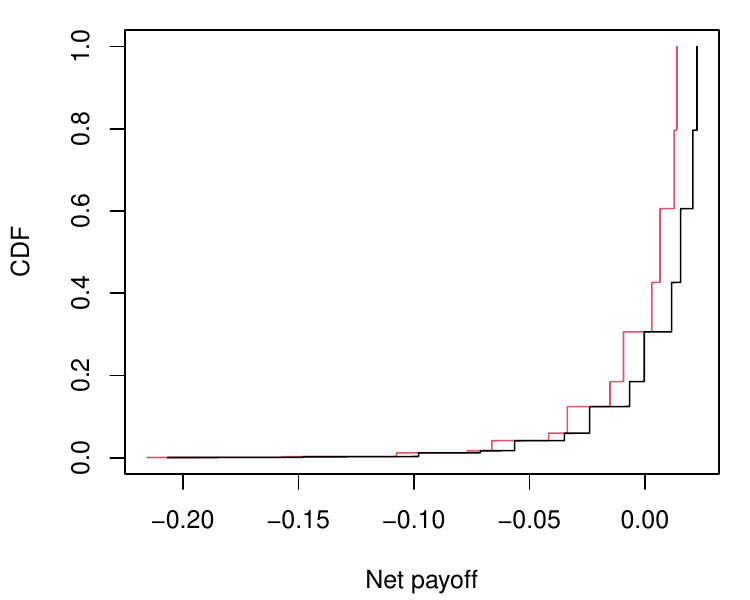}}
    \caption{Results for the binomial risk assessment comparison in Section~\ref{SS:bayes}.  Panel~(a) shows the IM (black), Bayes (red), and oracle (gray) risk assessments for $x=8$, $n=15$, and $\Theta=0.6$.  Panel~(b) shows the distribution of the IM (black) and Bayes (red) net payoffs, as a function of $X \sim \bin(n,\Theta)$ with $n=15$ and $\Theta=0.6$.}
    \label{fig:binom.risk}
\end{figure}

For a second concrete illustration, let $X \sim \bin(n,\Theta)$ as in Example~\ref{ex:binomial} above; let $\prior_x$ denote the Bayesian posterior based on the uniform prior for $\Theta$ and $\uPi_x$ the likelihood-based possibilistic IM.  For the comparison here, however, we opt for the simpler unscaled version of squared error loss, i.e., $\ell_a(\theta) = (a-\theta)^2$.  Figure~\ref{fig:binom.risk}(a) shows a plot of the Bayesian and IM risks when $n=15$ and $x=8$; also shown is the oracle assessment when $\Theta=0.6$.  While this only shows the case for a single $(n,x,\Theta)$ triplet, we hope the reader realizes that what's being depicted here is not uncommon.  Specifically, note how there's an entire interval of actions, roughly $a \in [0,0.4]$ for which the Bayes assessment is strictly less than the oracle assessment.  The IM assessment, however, never goes below the oracle assessment in this case.  To be clear, we {\em do not} claim that the IM assessment is literally always above the oracle assessment for every $(n,x,\Theta)$ pair, only that above-the-oracle is the general trend, as we demonstrate next.  We sampled 5000 data sets $X \sim \bin(n,\Theta)$ with $n=15$ and $\Theta=0.6$ and, for each, we evaluate the net payoffs $W_0^\text{\sc ba}(X; \Theta)$ and $W_0^\text{\sc im}(X; \Theta)$; Figure~\ref{fig:binom.risk}(b) shows the (empirical) distribution based on these 5000 samples.  What is immediately clear is that $W_0^\text{\sc im}(X; \Theta)$ is stochastically larger than $W_0^\text{\sc ba}(X; \Theta)$, which implies that the IMer's expected net payoff is greater than the Bayesian's.  What's not immediately clear is that the former expected value is positive ($0.0069$) whereas the latter is negative ($-0.0025$).  While the numerical values are small, the difference is strongly significant as verified by a simple t-test; we could make the gains/losses as large as we like since, in principle, changing the magnitude of the loss function doesn't affect acceptability/desirability.  The point is that, for some actions, the Bayesian risk assessment tends to exceed the oracle risk assessment and to the extent that the expected difference is negative, suggesting a systematic loss to the Bayesian making transactions that he deems to be acceptable (if not desirable).  If transactions that the Bayesian himself deems to be acceptable lead to systematic loss of money, then there's no other conclusion to draw but that the Bayesian assessment is overconfident.  


\subsection{Risk assessment properties: finite-sample}
\label{SS:finite}

As argued in the previous subsection, Bayesian and other precise probabilistic solutions tend to be overconfident.  That is, there are nontrivial sets of actions where the Bayes risk tends to be lower than the oracle risk.  Consequently, unbeknownst to the Bayesian, the price that he'd accept to take certain actions is systematically lower than what he'll lose by taking said action, hence a risk of systematically losing money.  This is a clear sign of unreliability.  The same numerical illustrations showing that the Bayes risk assessment tends to be overconfident for some actions also shows that the IM risk assessment doesn't tend to be overconfident in this sense.  That is, according to these (limited) empirical results, there are no sets of actions on which the IM's risk assessment is systematically lower than the oracle's.  In this section we prove that the IM's validity property, which implies reliability in inference, also implies a novel reliability property in decision-making: it's a rare event that there exists any action for which the IM's risk assessment is much smaller than an oracle's assessment. 
This, in turn, can be interpreted as a decision-theoretic version of the IM's validity property. Aside from offering a partial explanation for the numerical results in Section~\ref{SS:bayes} related to ``overconfidence,'' this is also suggestive of the overall quality of the IM's optimal action $\hat a(x)$ in \eqref{eq:best}: the IM's risk assessment effectively can't be smaller than the minimum oracle risk, so, if there was a set of actions on which $\uPi_x\ell_a$ tended to be too small compared to $\ell_a(\Theta)$, then it's likely that this gap would occur near $\hat a(x)$, where the IM risk is smallest, indicating that $\hat a(x)$ isn't near the oracle risk minimizer; but since the existence of such a set of actions is unlikely, it's safe to conclude that $\hat a(x)$ is not too far from the oracle risk minimizer.  

Toward comparing the IM and oracle assessments, a first basic question is if the two are even comparable.  Details will be provided in Section~\ref{SS:asymptotics} below, but the key point is that, under standard regularity conditions, the IM assessment $\uPi_{X^n} \ell_a$, depending on iid data $X=X^n$ of size $n$, converges to the oracle assessment $\ell_a(\Theta)$ in $\prob_\Theta$-probability as $n \to \infty$.  Therefore, it's not unreasonable to compare the IM and oracle assessments.

The general comparability between the two assessments suggests that we be more precise about what it means for the IM's assessment, $\uPi_x \ell_a$, to be ``too small'' compared to the oracle's assessment, $\ell_a(\Theta)$. Define the data-dependent local maximum loss, 
\begin{align}
L_a(x, \Theta) 
& = \sup\{ \ell_a(\theta): \pi_x(\theta) > \pi_x(\Theta)\}, \label{eq:max.loss.pl}
\end{align}
we'll refer to $L_a(x,\Theta)$ as the {\em quasi-oracle's assessment} of $a$, based on data $X=x$. This corresponds to an oracle who doesn't know $\Theta$ but knows which values $\theta$ are at least as $\pi_x$-plausible as $\Theta$ relative to data $X=x$. Then by ``too small'' we mean $\uPi_x \ell_a$ being less than a small multiple of $L_a(x,\Theta)$, i.e., 
\[ \uPi_x \ell_a \leq \alpha \, L_a(x,\Theta), \quad \alpha \in (0,1). \]
It's easy to see from \eqref{eq:max.loss.pl} that $L_a(x,\Theta) \geq \ell_a(\Theta)$, for all $x$, so the quasi-oracle's assessment is more conservative than the oracle's.  Therefore, if the IM's assessment, $\uPi_x \ell_a$, tends not to be too small compared to the quasi-oracle's assessment, $L_a(x,\Theta)$, then it also won't tend to be too small compared the oracle's assessment, $\ell_a(\Theta)$.  And not tending to be too small compared to the oracle suggests no overconfidence and provides some assurance that poor actions won't be favored by the IM.  The following theorem establishes the above decision-making reliability claim, namely, that it is indeed a rare event that the IM's assessment be too small compared to the quasi-oracle's.  A first version of this result appeared in \citet{imdec}, which was inspired by \citet{grunwald.safe} and, in turn, inspired similar results involving e-values in \citet{grunwald.epost}. 

\begin{theorem}
\label{thm:action}
Let $\ell_a: \TT \to [0,\infty)$ be a non-negative loss function, and define 
\begin{equation}
\label{eq:R.bound}
\ratio(x,\Theta) = \inf_{a \in \action} \frac{\uPi_x \ell_a - \ell_a(\hat\theta_x)}{L_a(x,\Theta) - \ell_a(\hat\theta_x)} \geq 0, 
\end{equation}
where $L_a(x,\Theta)$ is as in \eqref{eq:max.loss.pl} and $\hat\theta_x$ is the maximum likelihood estimator of $\Theta$.  Then 
\begin{equation}
\label{eq:R.bound0}
\sup_\theta \prob_{\theta}\{ \ratio(X,\theta) \leq \alpha \} \leq \alpha, \quad \text{for all $\alpha \in [0,1]$}. 
\end{equation}
This also implies the simpler-but-weaker result 
\[ \sup_\theta \prob_{\theta}\Bigl\{ \inf_{a \in \action} \frac{\uPi_X \ell_a}{L_a(X,\theta)} \leq \alpha \Bigr\} \leq \alpha, \quad \text{for all $\alpha \in [0,1]$}. \]
\end{theorem}

\begin{proof}
Define the function 
\begin{equation}
\label{eq:h.fun}
h_{x,a}(s) = \sup_{\theta: \pi_x(\theta) > s} \ell_a(\theta), \quad s \in [0,1], 
\end{equation}
so that the Choquet integral $\uPi_x \ell_a$ is just a Riemann integral of $h_{x,a}$.  It's clear that $s \mapsto h_{x,a}(s)$ is decreasing, which implies that, for any $\theta$, 
\begin{align*}
\uPi_x \ell_a & = \int_0^1 h_{x,a}(s) \, ds \\
& = \Bigl( \int_0^{\pi_x(\theta)} + \int_{\pi_x(\theta)}^1 \Bigr) \, h_{x,a}(s) \, ds \\
& \geq \pi_x(\theta) \, h_{x,a}(\pi_x(\theta)) + \{1 - \pi_x(\theta)\} \, h_{x,a}(1),
\end{align*}
where $h_{x,a}(1) = \ell_a(\hat\theta_x)$.  This is a Markov inequality for the Choquet integral, comparable to \citet{wang2011.choquet}.  Since $L_a(x,\theta) = h_{x,a}(\pi_x(\theta))$, it follows that
\begin{align*}
\ratio(X, \Theta) \leq \alpha & \iff \uPi_X \ell_a - \ell_a(\hat\theta_X) \leq \alpha \, \{L_a(X,\Theta) - \ell_a(\hat\theta_X)\} \quad \text{some $a$} \\
& \implies \pi_X(\Theta) \, \{ L_a(X,\Theta) - \ell_a(\hat\theta_X)\} \leq \alpha \, \{ L_a(X,\Theta) - \ell_a(\hat\theta_X) \}\quad \text{some $a$} \\
& \iff \pi_X(\Theta) \leq \alpha.
\end{align*}
Validity implies the latter event has $\prob_{\Theta}$-probability $\leq \alpha$, which proves \eqref{eq:R.bound}.  The simpler, second inequality \eqref{eq:R.bound0} holds because we can ignore the non-negative ``$\{1-\pi_x(\theta)\} h_{x,a}(1)$'' in the lower bound for $\uPi_x\ell_a$ stated above.  
\end{proof}

The uniformity in the action $a$ implies that the same result holds with a data-driven action $a(x)$ plugged in, as the following corollary shows.  This includes plugging in $\hat a(x) = \arg\min_a \uPi_x \ell_a$, the IM risk minimizer in \eqref{eq:best}.  

\begin{corollary}
\label{cor:action.min}
Under the setup of Theorem~\ref{thm:action}, if $a: \XX \to \action$ is a data-dependent action, e.g., the IM risk minimizer in \eqref{eq:best}, then 
\[ \sup_\theta \prob_{\theta}\{ \uPi_X \ell_{a(X)} \leq \alpha \, L_{a(X)}(X,\theta) \} \leq \alpha, \quad \alpha \in [0,1]. \]
\end{corollary}

\begin{proof}
The infimum ratio over actions $a \in \action$ can be no larger than the ratio at a specific (data-dependent) action $a(x)$.  Then the claim follows from \eqref{eq:R.bound0}.  
\end{proof}

The ``quasi-oracle risk assessment'' complicates interpretation but, fortunately, it's easy to remove this complication by replacing the quasi-oracle assessment $L_a(X,\Theta)$ with the simpler oracle assessment $\ell_a(\Theta)$.  

\begin{corollary}
\label{cor:action.oracle}
Under the setup of Theorem~\ref{thm:action}, 
\begin{equation}
\label{eq:action.oracle.simplest}
\sup_\theta \prob_\theta\Bigl\{ \inf_{a \in \action} \frac{\uPi_X \ell_a}{\ell_a(\theta)} \leq \alpha \Bigr\} \leq \alpha, \quad \alpha \in [0,1]. 
\end{equation}
Furthermore, just as in Corollary~\ref{cor:action.min}, if $a(X)$ is a data-driven action, then 
\[ \sup_\theta \prob_{\theta}\{ \uPi_X \ell_{a(X)} \leq \alpha \, \ell_{a(X)}(\theta) \} \leq \alpha, \quad \alpha \in [0,1]. \]
\end{corollary}

\begin{proof}
This follows from Theorem~\ref{thm:action} and the fact that $\ell_a(\theta) \leq L_a(x,\theta)$ for all $(x,\theta,a)$. 
\end{proof}

It's not worth stating as another formal corollary, but there's one other potentially helpful way to express the above results---in terms of expected loss differences rather than ratios.  Let's take the simplest case as in Corollary~\ref{cor:action.oracle}.  Then 
\[ \frac{\uPi_x \ell_a}{\ell_a(\theta)} \leq \alpha \iff \uPi_x\ell_a - \ell_a(\theta) \leq -(1-\alpha) \ell_a(\theta), \]
and, therefore, \eqref{eq:action.oracle.simplest} is equivalent to 
\[ \sup_\theta \prob_\theta \bigl\{ \uPi_X\ell_a - \ell_a(\theta) \leq -(1-\alpha) \ell_a(\theta) \; \text{ for some $a \in \action$} \bigr\} \leq \alpha. \]
This equivalent event corresponds to the difference between what the IMer user in exchange for taking action $a$ and what she actually loses when $\Theta=\theta$ being strictly negative.  It's a probability similar to the one above that was plotted in Figure~\ref{fig:binom.risk} for the binomial example.  So, the result in the above display helps to explain why the IM's assessment is safe from the kind of ``overconfidence'' that the Bayesian assessment is afflicted by. 

For a quick recap, recall that it would be undesirable if the IM's assessment of $a$ were too small compared to the oracle's or the quasi-oracle's since it would put the decision-maker at risk of suffering non-trivial loss.  This is especially true for the ``best'' action, $a=\hat a(x)$, the one that the decision-maker is likely to take.  By Theorem~\ref{thm:action}, the IM's validity implies that such undesirable cases are rare events with respect to the  distribution of $X$.  Therefore, validity provides some assurance that the IM's data-driven assessment of action $a$ is not inconsistent with that of the oracle or quasi-oracle, hence the IM helps the decision-maker mitigate risk. 



\subsection{Risk assessment properties: large-sample}
\label{SS:asymptotics}

Section~\ref{SS:bayes} showed, mostly empirically, that Bayes risk assessments tend to be too small and hence overconfident compared to the oracle assessment.  The previous section demonstrated theoretically that the IM risk assessment is larger and hence not prone to the same form of overconfidence; in fact, there's even a certain calibration between the IM's and quasi-oracle's risk.  Of course, it would be easy to avoid overconfidence by going too far in the other direction, by being grossly conservative.  The result in Theorem~\ref{thm:limit} below shows that this is not how the IM solution achieves its reliability.  In fact, just like the Bayes risk, the IM risk merges with the oracle risk asymptotically, as the sample size $n$ approaches infinity, so the differences with respect to overconfidence can't be explained by the IM risk being significantly too large/conservative.  

The following large-sample result builds on the recent developments in \citet{imbvm.ext}.  Like all general large-sample results, certain regularity conditions must be imposed on the model $\{\prob_{\theta}: \theta \in \TT\}$ to ensure that troublesome anomalies and edge cases are avoided.  The standard set of regularity conditions are commonly referred to as the {\em Cram{\'e}r conditions} \citep{cramer.book}, versions of which can be found in classical texts, including \citet[][Theorem~3.10]{lehmann.casella.1998} and \citet[][Theorem~7.63]{schervish1995}. More modern treatments \citep[e.g.,][]{vaart1998, bickel1998} adopt the weaker conditions of \citet{lecam1970}, as we do here. Roughly, these conditions ensure that the density/mass functions $p_\theta$ corresponding to the $\prob_\theta$s are sufficiently smooth so that $\theta \mapsto \log p_\theta(x)$ admits a suitable quadratic approximation; it's from this quadratic approximation of the log-likelihood that asymptotic Gaussianity is derived. Details are provided in Section~3.2 of \citet{imbvm.ext}, and a summary is given in Appendix~\ref{AA:limit} below. In particular, it's under precisely these conditions that \citet{imbvm.ext} established a large-sample Gaussian approximation of the likelihood-based possibilistic IM, a possibility-theoretic version of the celebrated {\em Bernstein--von Mises theorem} in Bayesian statistics \citep[e.g.,][Theorem~10.1]{vaart1998}. 

Theorem~\ref{thm:limit} that follows imposes these standard regularity conditions along with some very basic restrictions on the loss function.  This is definitely not the best result possible, in terms of both the strength of the assumptions and the breadth of application; for example, the theorem doesn't apply to testing with zero--one loss, in part because this would be fairly easy to handle directly.  Our moderate goal here is to present a single result that is both relatively simple and general, while relying on the efforts in \citet{imbvm.ext} for most of the heavy lifting.  


\begin{theorem}
\label{thm:limit}
Assume the regularity conditions as in Section~3.2 of \citet{imbvm.ext} and summarized in Appendix~\ref{AA:limit} below.  For each $a \in \action$, let $\ell_{a}$ be a non-negative and continuous function with compact support.   Then $\uPi_{X^n} \ell_a \to \ell_{a}(\Theta)$ in $\prob_{\Theta}$-probability as $n \to \infty$, for each action $a \in \action$. Moreover, if the collection $\{\ell_a: a \in \action\}$ is uniformly bounded and equicontinuous, then the convergence is uniform in actions, i.e., $\sup_{a \in \action} \bigl| \uPi_{X^n} \ell_a - \ell_a(\Theta) \bigr| \to 0$ in $\prob_\Theta$-probability as $n \to \infty$. 
\end{theorem}

\begin{proof}
See Appendix~\ref{AA:limit}. 
\end{proof}

Again, some discussion of the aforementioned regularity conditions is given in Appendix~\ref{AA:limit} below.  The ``compact support'' condition is an obstacle in many practical situations, but this definitely isn't a necessary condition in principle; we need it here only because the IM's Bernstein--von Mises theorem comes with only a ``compact convergence'' or ``converges uniformly on compact sets'' conclusion.  Finally, equicontinuity is implied by some more basic and familiar conditions, e.g., if $\theta \mapsto \ell_a(\theta)$ is Lipschitz continous with Lipschitz constant bounded in $a$.  

To summarize, asymptotic merging of the IM risk $\uPi_{X^n} \ell_a$ with the oracle risk $\ell_a(\Theta)$ implies that there is no grossly conservative cushion preventing the IM from being overconfident in the same way that the Bayes risk is.  In fact, under the same conditions, the Bayes risk $\prior_{X^n} \ell_a$ also merges with the oracle, so all three will agree in the idealistic, perfect-information, large-sample limit.  In the practically realistic cases with imperfect information, however, only the IM is reliable in the sense of Theorem~\ref{thm:action}.

\section{Applications}
\label{S:applications}

\subsection{Gamma model}
\label{SS:gamma}

\begin{figure}[t]
    \centering
    \subfigure[Contour function]{\includegraphics[width=0.49\linewidth]{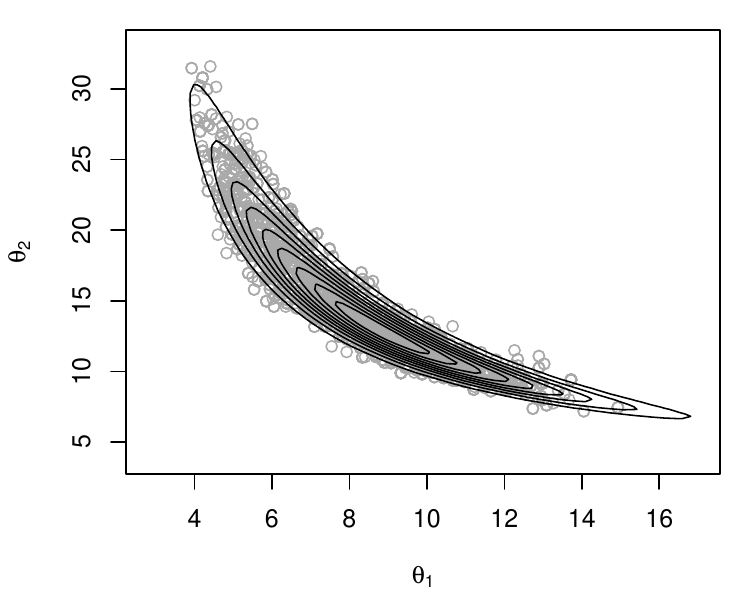}}
    \subfigure[Risk assessments]{\includegraphics[width=0.49\linewidth]{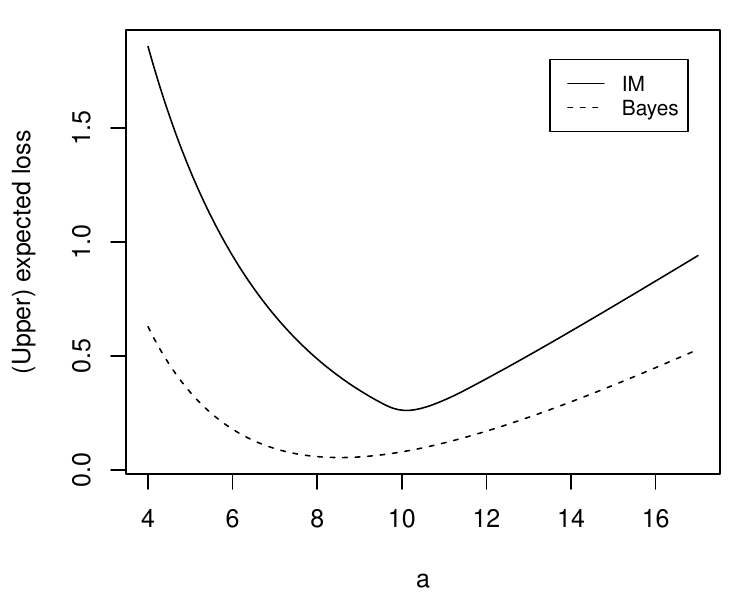}}
    \caption{Results from the gamma illustration in Section~\ref{SS:gamma}.  Panel~(a) compares the Jeffreys-prior Bayes posterior samples (gray) with the possibilistic IM contour (black).  Panel~(b) shows the IM (solid) and Bayes (dashed) risk assessments.}
    \label{fig:gamma}
\end{figure}

Let $X=(X_1,\ldots,X_n)$ consists    of independent and identically distributed random variables from a gamma distribution with unknown shape parameter $\Theta_1 > 0$ and unknown scale parameter $\Theta_2 > 0$.  The gamma distribution is commonly used as a model for time-to-event data, e.g., in survival or reliability analysis.  Neither the Jeffreys-prior Bayesian solution nor the possibilistic IM solution have closed-form expressions, but summaries can easily be found numerically.  For illustration, consider the data in Example~3 of \cite{fraser.reid.wong.1997} on the survival time (in weeks) for $n=20$ rats exposed to a certain level of radiation.  
Figure~\ref{fig:gamma}(a) shows a sample from the Bayesian posterior distribution of $\Theta = (\Theta_1, \Theta_2)$, given $X=x$, with the possibilistic IM contour overlaid.  For the decision-making problem, we focus here on estimating the shape using the loss 
\[
\ell_a(\theta) = ( \log a - \log \theta_1 )^2, \quad a > 0.
\]
Figure~\ref{fig:gamma}(b) plots the Bayesian and IM risk assessments as a function of the action $a$.  Naturally, since the two ``posterior distributions'' are visually similar, the shapes of the two risk assessments are also similar.  But the two differ in both the magnitudes and the location of the minimizer.  That the Bayes risk assessment is quite small compared to that of the IM across the entire range of actions is consistent with the ``Bayes is overconfident'' claim offered in Section~\ref{SS:bayes}.  

\subsection{Random effects model}
\label{SS:vc}

Random effect models are ubiquitous in applications across business, engineering, and science.  The simplest such model is commonly written as 
\[ X_{ij} = \theta_0 + T_i + E_{ij}, \quad i=1,\ldots,m, \quad j=1,\ldots,n, \]
where $X=(X_{ij})$ is the observable data and $T_i$ and $E_{ij}$ are mutually independent random effects, with $E_{ij} \iid \nm(0, \theta_1)$ and $T_i \iid \nm(0,\theta_2)$.  Then the three-dimensional parameter that indexes this model is $\theta=(\theta_0, \theta_1, \theta_2)$, where $\theta_0$ is the overall mean and $(\theta_1,\theta_2)$ are variance components associated with error/replication and treatment, respectively.  The overall mean is not of primary interest, and it can be easily marginalized out, so the focus will be on the variance components, in particular, that for the treatment effect. 

\begin{figure}[t]
    \centering
    
    \subfigure[Contour function]{\includegraphics[width=0.49\linewidth]{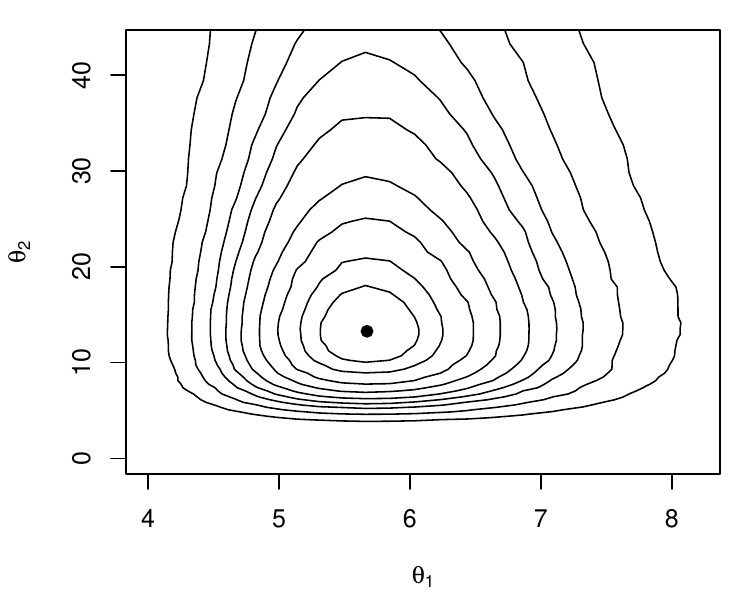}}
    \subfigure[Risk assessments]{\includegraphics[width=0.49\linewidth]{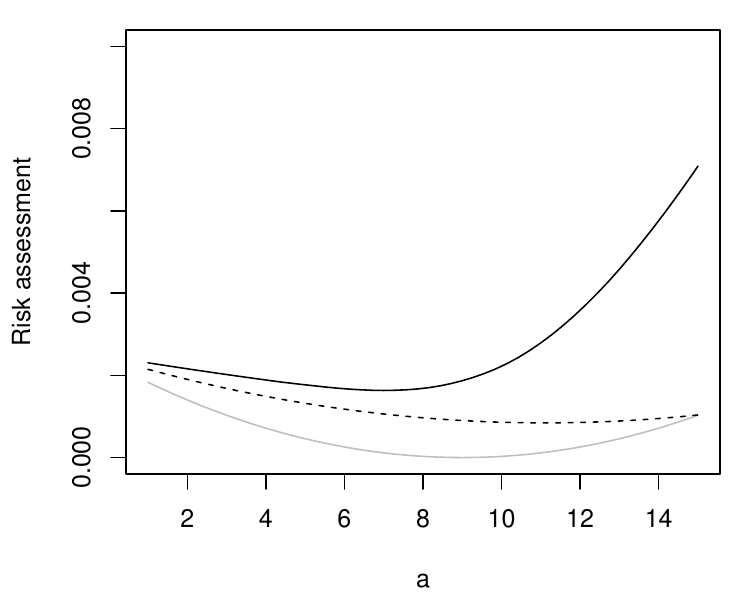}}
    \caption{Results from the random effects model in Section~\ref{SS:vc}. Panel (a) shows the possibilistic IM contour for the variance components $(\theta_1, \theta_2)$. Panel (b) compares the oracle (gray), Bayes (dashed), and IM (solid) risk based on loss function \eqref{eq:vc.loss}.}
    \label{fig:vc}
\end{figure}

The work that follows is based on a simulated data set with $m=5$, $n=20$, and true parameters $(\Theta_0, \Theta_1, \Theta_2) = (0, 7, 9)$.
A plot of the contour $\pi_x$, as a function of the variance components $(\theta_1,\theta_2)$, is shown in Figure~\ref{fig:vc}(a). The dot in the center marks the mode of $\pi_x$, which is the maximum likelihood estimator, $\hat\theta = (5.67, 13.28)$.  

For estimation of the treatment variance component, \citet{portnoy1971} suggests the following variation on the usual squared-error loss function:
\begin{equation}
\label{eq:vc.loss}
\ell_a(\theta) = \frac{(\theta_2 - a)^2}{(\theta_1 + n \theta_2)^2}, \quad a > 0. 
\end{equation}
There are different ways to assess the risk associated with actions relative to this loss function.  One is the oracle assessment, another is 
the Jeffreys-prior Bayes assessment as in \citep{tiao.tan.1965}, and the IM assessment proposed here.  A plot of those assessments is shown in Figure~\ref{fig:vc}(b).  The key observation is that these all have similar shapes and magnitudes, as expected; the IM assessment favors an action smaller than the actions favored by the other assessments, which is consistent with both the loss function and IM solution's inherent reliability/safety properties. 
Further, to compare the Bayes and IM estimators, we ran multiple experiments and found that the frequentist expected loss of the Bayes and of the IM estimators are equivalent. 

\subsection{Logistic regression model}
\label{SS:logistic}

Thus far we have focused on continuous loss functions and different variations on the squared loss.  But recall that the hypothesis testing scenario in Section~\ref{SS:bayes}  with 0--1 loss revealed non-negligible differences between the Bayes and IM risks, as clearly shown in Figure~\ref{fig:gauss.test}. The logistic regression model discussed here reveals similar differences.  To set the scene, for each $i=1,\ldots,n$, the model $\prob_\theta$ posits that  
\[
(Y_i \mid X_i=x_i) \sim {\sf Ber}\{ f_\theta(x_i) \} \quad \text{where} \quad f_\theta(x) = \bigl(1 + e^{-(\theta_0 + \theta_1 x)} \bigr)^{-1}. 
\]
The standard Bayesian solution introduces a prior for the coefficients $\Theta=(\Theta_0, \Theta_1)$; for the present example, we use an exchangeable prior for $\Theta$ defined as a mean-zero Gaussian scale mixture as suggested in, e.g., \citet{gelman.weak.2008}; see, also, \citet[][Sec.~19]{yang.berger.1998} for other default priors.  Neither the Bayes nor IM quantification of uncertainty has a closed-form expression, so Monte Carlo methods are needed for both.  

The decision problem here is related to prediction/classification.  Suppose there's a fixed value $x$ of the covariate and the action $a \in \{0,1\}$ is a guess of the value of $y$ associated with the given $x$.  A common strategy is, for given $x$ and $\theta$, to assign category $a=0$ if $f_\theta(x) < \frac12$ and category $a=1$ if $f_\theta(x) \geq \frac12$.  Note that the two conditions reduce to $\theta_0 + \theta_1 x < 0$ and $\theta_0 + \theta_1 x \geq 0$.  Then the loss function that assigns a common cost of 1 unit for misclassification in either direction, i.e., 
\[ \ell_a(\theta) = 1(\theta_0 + \theta_1 x \geq 0) \, 1(a=0) + 1(\theta_0 + \theta_1 x < 0) \, 1(a=1). \]
In what follows, we compare the Bayes and IM expected losses, for each $a \in \{0,1\}$, treated as functions of the underlying covariate value $x$.  

\begin{figure}[t]
    \centering
    
    \subfigure[Chloroacetic acid data]{\includegraphics[width=0.49\linewidth]{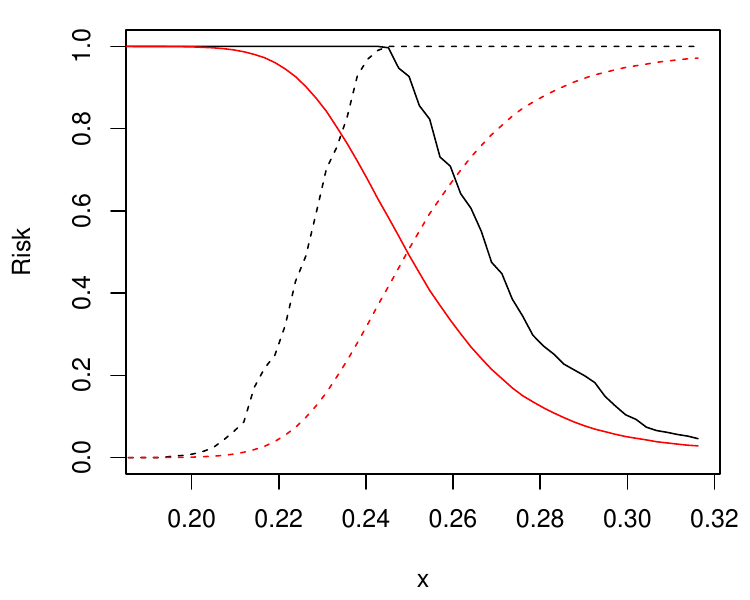}}
    \subfigure[Simulated data]{\includegraphics[width=0.49\linewidth]{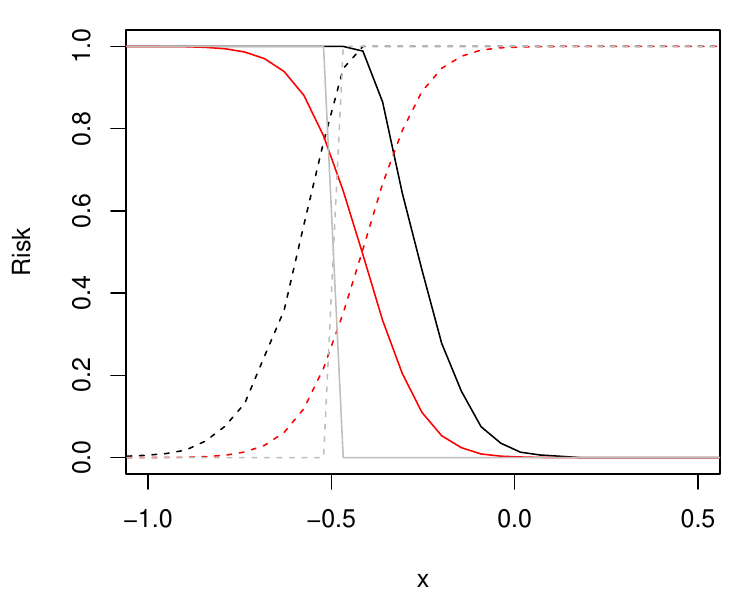}}
\caption{Results from the logistic regression model in Section~\ref{SS:logistic}. Panel (a) compares risk assessment from Bayes (red) and IM (black) constructions for data about the relationship between chloroacetic acid exposure and mouse mortality. Panel (b) compares risk assessment from oracle (gray), Bayes (red), and IM (black) in simulated data. The solid lines denote when $a=1$ and dashed lines are for $a=0$.}
    \label{fig:loss_logistic}
\end{figure}

We first consider data from Table 8.4 of \cite{ghosh-etal-book}, which concerns the relationship between the exposure to chloroacetic acid and mouse mortality. Figure~\ref{fig:loss_logistic}(a) shows that the IM risk assessment is strictly larger than the Bayes assessment for both $a=0$ and $a=1$, which is desirable behavior in light of Bayes tending to be overconfident as discussed in Section~\ref{SS:bayes}.  To compare risk assessments to the oracle, we set $(\Theta_0, \Theta_1) = (1,2)$ to simulate data with $n=100$ and the $X_i$s iid $\nm(0, 1)$.  Figure~\ref{fig:loss_logistic}(b) demonstrates that the IM risk is again strictly greater than the Bayes risk and, moreover, consistently agrees with the oracle. 
This result reiterates the finding of Section~\ref{SS:bayes}, namely, that the overconfident Bayesian tends to under-estimate the oracle risk, making him susceptible to losing money, whereas the IMer has apparently just the right amount of caution to protect herself from systematically falling short of the oracle. 

\section{Invariant decision problems}
\label{S:group}

\subsection{Motivation}

The previous sections focused on showing that the IMer's assessments of different actions are not overconfident in the same way that Bayesian's assessments tend to be overconfident, as demonstrated in Section~\ref{SS:bayes}.  This focus on being ``not overconfident,'' together with our recommendation for use of a minimax strategy, suggests that the IM solution is perhaps too conservative.  This is not true, however, as we've already highlighted in some specific examples above: indeed, the IM solutions do sometimes agree with the Bayesian solutions, which in turn are sometimes optimal with respect to various classical decision-theoretic considerations, e.g., minimum variance unbiased estimators or minimum risk equivariant estimators.  The goal of the present section is to formalize the results in those examples, to offer sufficient conditions on when the IM solution will, in fact, agree with what is generally considered to be the optimal solution (which happens to be a certain Bayesian solution in this context).  So, the reader can interpret the results in this section as an efficiency-focused complement to those in Section~\ref{SS:finite} above establishing that the IM solution is always safe and reliable, and a finite-sample complement to the large-sample efficiency result in Section~\ref{SS:asymptotics}.


\subsection{Setup}

Let $\G$ denote a group of bijections $g: \XX \to \XX$ acting on the sample space $\XX$, with function composition $\circ$ as the binary operation. 
As is customary in the literature, We'll write $gx$ for the image of $x \in \XX$ under transformation $g \in \G$.  Since $\G$ is a group, it's associative, i.e., $g_1 \circ (g_2 \circ g_3) = (g_1 \circ g_2) \circ g_3$ for all $g_1,g_2,g_3 \in \G$, it contains the identity transformation, and for every $g \in \G$, there exists an inverse $g^{-1} \in \G$ such that $g \circ g^{-1} = g^{-1} \circ g = \text{identity}$. Examples include translations, scale changes, rotations, and permutations.  

The group $\G$ connects to the statistical model as follows.  Suppose that, for each $g \in \G$ and each $\theta \in \TT$, there exists a corresponding $\bar g \theta \in \TT$ such that 
\begin{equation}
\label{eq:invariant}
\prob_\theta(gX \in \cdot) = \prob_{\bar g \theta}(X \in \cdot), \quad (\theta,g) \in \TT \times \G. 
\end{equation}
For example, if the distribution of $X$ depends on a location parameter $\theta$, then the distribution of $X+\tau$ depends on parameter $\theta + \tau$.  When the statistical model $\{\prob_{\theta}: \theta \in \TT\}$ satisfies \eqref{eq:invariant}, we follow \citet[][Ch.~6]{berger1985}, \citet[][Ch.~3]{eaton1989}, and \citet[][Ch.~6]{schervish1995}, and say that it's an {\em invariant statistical model} in the sense that the model itself remains unchanged when data are transformed by elements of the group.  It's often the case in applications that the distributions $\prob_{\theta}$ have densities respect to some underlying dominating measure, and then invariance means 
\begin{equation}
\label{eq:invariant.density}
p_\theta(x) = p_{\bar g \theta}(gx) \, \chi(g), \quad x \in \XX, \; \theta \in \TT, \; g \in \G, 
\end{equation}
where $\chi(g)$ is the ``multiplier,'' a change-of-variables Jacobian term.  For details, see \citet[][Ch.~6]{schervish1995}, \citet{eaton1989}, and many others.  

Define $\Gbar$ as the collection of all those $\bar g: \TT \to \TT$, corresponding to the mappings $g \in \G$; this too is a group.  A further simplification, commonly adopted in the literature \citep[e.g.,][p.~371]{schervish1995}, is to assume $\TT = \G = \Gbar$. This assumption simply means we don't have to distinguish $g$, $\bar g$, and $\theta$ or define functions that connect them.  

The set $\G x = \{gx: g \in \G\} \subseteq \XX$ is called the {\em orbit} of $\G$ corresponding to $x$.  The orbits partition $\XX$ into equivalence classes, so every point $x \in \XX$ falls on exactly one orbit.  This partition can be used to construct a new coordinate system on $\XX$ which will be useful for us in what follows.  Identify $x \in \XX$ with $(g_x, u_x)$, where $u_x \in \UU$ denotes the label of orbit $\G x$ and $g_x \in \G$ denotes the position 
of $x$ on the orbit $\G x$.  

One of the key results in this section is a certain correspondence between the ``standard'' Bayesian (and fiducial) solution and the possibilistic IM solution in this invariant model context.  For this correspondence to make any sense, the group with respect to which the model is invariant must be unique.  Otherwise, if there are multiple groups with respect to which the model is invariant, then each group induces its own invariant prior and, in turn, its own corresponding Bayes solution.  Since it's not possible for the IM solution to agree with several different Bayesian solutions simultaneously, we assume throughout that the group in question is unique.  Note that uniqueness is not universal: for example, in covariance matrix estimation \citep[e.g.,][]{yang.berger.1994} there are different groups with respect to which the model is equivariant: the standard choice is left-multiplication by lower triangular matrices with positive diagonal entries, but invariance also holds with respect to multiplication by arbitrary positive definite matrices.  But keep in mind that uniqueness of the group, which is a rather strong condition, is not a constraint on the IM solution or its properties; the uniqueness is needed only to form the connection between existing Bayesian (and fiducial) solutions.  In addition to uniqueness of the group, we make the following 

\begin{asmp}
Let $\{p_\theta: \theta \in \TT\}$ be a family of densities invariant with respect to a locally compact topological group $\G$ in the sense of \eqref{eq:invariant.density} and, as explained above, take $\TT = \G = \Gbar$.  In addition, the following will be assumed:
\begin{itemize}
\item[A1.] The left Haar measure $\lambda$ and the corresponding right Haar measure $\rho$ on (the Borel $\sigma$-algebra of) $\G$ exist and are unique up to scalar multiples.
\item[A2.] There exists a bijection $t: \XX \to \G \times \UU$, with both $t$ and $t^{-1}$ measurable, that maps $x \in \XX$ to its position--orbit coordinates $(g_x, u_x) \in \G \times \UU$.  
\item[A3.] The distribution of $t(X) = (G, U) \in \G \times \UU$ induced by the distribution of $X \sim \prob_{\theta}$ has a density with respect to $\lambda \times \mu$ for some measure $\mu$ on $\UU$.
\end{itemize} 
\end{asmp}

A few quick remarks are in order.  First, for A1, existence and uniqueness of the left and associated right Haar measures on locally compact topological groups is a classical result \citep[e.g.,][]{halmos.measure, nachbin1965}.  For A2, note that $t(gx) = (g \circ g_x, u_x)$ for all $g \in \G$.  That is, $g$ only acts on the first coordinate in $t$, so it's invariant with respect to $\G$ in the second coordinate.  For A3, existence of a joint density with respect to a product measure ensures that there is no difficulty in defining a conditional distribution for $G$, given $U=u$.  Finally, $U=U_X$, as a function of $X$, is an ancillary statistic, i.e., its distribution under $\prob_\theta$ doesn't depend on $\theta$. 

The simplest example is a   location model where $\G = \Gbar = (\RR, +)$.  This is an Abelian group, so the left and right Haar measures are the same and both equal to Lebesgue measure.  The function $x \mapsto t(x)$ in A2 consists of two components: in its ``$g_x$'' coordinate an equivariant function of $x$ that estimates the location and, in its ``$u_x$'' component, an invariant function of $x$, such as residuals.  For example, $g_x = \bar x$ the arithmetic mean of $x=(x_1,\ldots,x_n)$ and $u_x = \{x_i - \bar x: i=1,\ldots,n\}$.  Note that the $u_x$ coordinate satisfies a constraint, so, after it's represented in a suitable lower-dimensional space $\UU$, the measure $\mu$ can be taken as Lebesgue measure there.  Many other problems fit this general form \citep[e.g.,][Ch.~1--2, including exercises]{fraser1968}. 

Finally, when a suitable ancillary statistic can be identified, like $U$ above, the recommendation in \citet{martin.partial2}---akin to recommendations by \citet{fisher1973}, \citet{berger.copss}, and others---is to reduce the complexity/dimension by conditioning on the observed value of the ancillary statistic in the validification step \eqref{eq:contour}.  We will adopt this conditional version of the IM construction below; see \eqref{eq:cond.im}.


\subsection{Bayes--fiducial--IM connection}
\label{SS:connection}

When the statistical model is invariant with respect to a group of transformations, an interesting connection between the familiar Bayesian/fiducial solution and the IM solution emerges.  This connection is through the IM's credal set.  

Recall that we're assuming $\TT = \G = \Gbar$.  In this case, generic values of $\theta \in \TT$ can be identified with transformations in $\G$; the same goes for the uncertain variable $\Theta$.  So, in what follows, we treat $\theta$ as a transformation on $\XX$ that can be inverted to $\theta^{-1}$ and can be composed via $\circ$ with other transformations in $\G$.  A key result \citep[][Cor.~6.64]{schervish1995} is that the density of $X$ or, equivalently, of $t(X) = (G,U)$, under $\prob_{\theta}$, is given by 
\begin{equation}
\label{eq:GU.joint}
p_\theta(g, u) = f(\theta^{-1} \circ g, u), 
\end{equation}
where the function $f: \G \times \UU \to \RR$ doesn't directly depend on $\theta$.  The particular form of $f$ isn't important---all that matters is how the right-hand side above depends on $\theta$.  Next are two important consequences of \eqref{eq:GU.joint}.
\begin{itemize}
\item Lemma~6.65 in \citet{schervish1995} says that the Bayesian posterior distribution for $\Theta$, given $x \equiv (g,u)$, has a density with respect to the right Haar prior measure 
\begin{equation}
\label{eq:schervish}
q_x^\star(\theta) = c_g \, c_u' \, f(\theta^{-1} \circ g, u), 
\end{equation}
where $c_g$ and $c_u'$ only depend on the $g$ and $u$-components of $x$, respectively.  This agrees with the accepted fiducial distribution for $\Theta$, given $x \equiv (g,u)$, for equivariant models like assumed here. This generalizes the classical result of \citet{lindley1958} on the connection between Bayes and fiducial; this is also the context around which Fraser's ``structural inference'' is developed \citep{fraser1968, fraser1965, fraser1966a}. 
\item Proposition~2 in \citet{martin.isipta2023} shows that the relative likelihood is given by $R(x,\theta) = d_u \, f(\theta^{-1} \circ g, u)$, where $d_u$ depends only on the $u$-component of $x$.  
\end{itemize} 
These points are importance for the following calculation.  Write $\pi_{g|u}$ for the possibilistic IM's contour, where the subscript emphasizes that the conditional distribution of $X \equiv (G,U)$, given the observed value of $u$, is used in the validification step \eqref{eq:contour}.  Then 
\begin{align}
\pi_{g|u}(\theta) & := \prob_\theta\{ R(X,\theta) \leq R(x,\theta) \mid U=u \} \label{eq:cond.im} \\
& = \prob_\theta\{ f(\theta^{-1} \circ G, u) \leq f(\theta^{-1} \circ g, u) \mid U=u \} \notag \\
& = \prob\{ f(H, u) \leq f(\theta^{-1} \circ g, u) \mid U=u \} \notag \\
& = \prior_x^\star\{ f(H, u) \leq f(\theta^{-1} \circ g, u) \} \notag \\
& = \prior_x^\star\{ q_x^\star(\Theta) \leq q_x^\star(\theta) \}, \notag
\end{align}
where the first line is by definition; the second line is by simplification; the third line is by the fact that $H := \theta^{-1} \circ G$ is a pivot with respect to the conditional distribution of $G$, given $U=u$, under $\prob_\theta$; the fourth line is by the fact that $H := \Theta^{-1} \circ g$ is also a pivot with respect to the Bayes posterior $\Theta \sim \prior_x^\star$ and has the same distribution as $H$ in previous line; and the last line is by simplification via \eqref{eq:schervish}. 

The last line in the above display can be recognized as the probability-to-possibility transform applied to the Bayesian/fiducial posterior distribution.  This observation leads to a new proof of the main result in \citet{martin.isipta2023}.

\begin{theorem}
\label{thm:fiducial}
Under the equivariant model setup described above, the Bayesian/fiducial posterior distribution $\prior_x^\star$ is the maximal element in the credal set $\cred(\uPi_{g|u})$ associated with the possibilistic IM.  That is, $\prior_x^\star$ satisfies 
\begin{equation}
\label{eq:Qstar.equal}
\prior_x^\star\bigl( \{\theta: \pi_{g|u}(\theta) > \alpha\} \bigr) = 1-\alpha, \quad \alpha \in [0,1]. 
\end{equation}
\end{theorem}

The take-away message is that, under an equivariant model, the familiar Bayes/fiducial distribution has a new interpretation as the maximally diffuse element in the possibilistic IM's credal set---what \citet{immc} calls the inner probabilistic approximation of the possibilistic IM.  For example, this connection explains why the Bayes/fiducial solutions offer exact confidence regions in equivariant models but not more generally: IMs are always valid and, hence, yield exact confidence regions, but the familiar probabilistic solutions can only do so when they match up close enough to the IM solution, which they do in equivariant problems.  See, also, \citet{reimagined}.

\subsection{Risk assessments}
\label{SS:invariant.risk}

The loss function is invariant if for every $\bar g \in \Gbar$ there exists $\tilde g: \action \to \action$ with
\[ \ell_{\tilde g a}(\bar g \theta) = \ell_a(\theta) \quad \text{all $(a,\theta)$}. \]
The collection $\Gtilde$ of all such maps $\tilde g$ also forms a group under function composition; just as above, for simplicity, we'll assume below that $\G = \Gbar = \Gtilde = \TT$.  

A well-known fact \citep[e.g.,][Theorem~6.59]{schervish1995} is that, when the loss function is invariant, the minimum risk equivariant decision rule is the formal Bayes rule, the minimizer of the expected loss with respect to the posterior distribution based on the right Haar prior $\rho$.  In the present notation, this optimal decision rule is 
\begin{equation}
\label{eq:bayes.rule}
\hat a(x) = \inf_{a \in \action} \prior_x^\star \ell_a. 
\end{equation}
We show below that, under an additional condition on the loss function, the same $\hat a(x)$ minimizes the IM's upper expected loss and, moreover, the IM and Bayes risk assessments of $\hat a$ are identical, i.e., $\uPi_x \ell_{\hat a(x)} = \prior_x^\star \ell_{\hat a(x)}$.  Here we're returning to our previous and simpler notation $\uPi_x$, but note that this stands for $\uPi_{g|u}$ as in Section~\ref{SS:connection}.  

The specific condition we assume about the loss function is unfamiliar, but not uncommon.  Define a data-dependent partial order, $\preceq_x$, on $\TT$ as follows:
\begin{equation}
\label{eq:order}
\theta \preceq_x \theta' \iff \pi_x(\theta) \leq \pi_x(\theta'). 
\end{equation}
Then the maximum likelihood estimator is the ``$\preceq_x$-largest'' element of $\TT$ and the more confidence sets $\{C_\alpha(x): \alpha \in [0,1]\}$ in \eqref{eq:region} that an element $\theta$ belongs to, the ``$\preceq_x$-larger'' it is.  Then the specific condition we impose is that, for $\hat a(x)$ as defined in \eqref{eq:bayes.rule}, the loss function $\theta \mapsto \ell_{\hat a(x)}(\theta)$ is $\preceq_x$-increasing, i.e., 
\begin{equation}
\label{eq:increasing}
\theta \preceq_x \theta' \implies \ell_{\hat a(x)}(\theta) \leq \ell_{\hat a(x)}(\theta'). 
\end{equation}
The simplest example of such a loss function is squared error loss in a symmetric location parameter model.  A general investigation into the properties of probability distribution families under partial orders like the one above can be found in \citet{bergmann1991}.  

\begin{theorem}
\label{thm:risk}
Consider an equivariant statistical model, with $\uPi_x$ the IM's upper probability and the Bayes posterior $\prior_x^\star$ based on the right Haar prior the maximal element in the IM's credal set.  If the loss function is invariant and, for the Bayes rule $\hat a(x)$ in \eqref{eq:bayes.rule}, satisfies \eqref{eq:increasing} relative to the partial ordering \eqref{eq:order}, then 
\begin{equation}
\label{eq:risk}
\inf_{a \in \action} \uPi_x \ell_a = \uPi_x \ell_{\hat a(x)} = \prior_x^\star \ell_{\hat a(x)}. 
\end{equation}
Therefore, the Bayes rule $\hat a(x)$ is also the IM's minimum upper expected loss action, and the corresponding risk assessments at $\hat a(x)$ are the same.  
\end{theorem}

\begin{proof}
By $\preceq_x$-monotonicity in \eqref{eq:increasing}, there exists a function $\beta_x(\cdot)$, whose form isn't important, such that 
\begin{equation}
\label{eq:increasing2}
\ell_{\hat a(x)}(\theta) > t \iff \pi_x(\theta) < \beta_x(t), \quad t \geq 0. 
\end{equation}
Using the basic definition of the Choquet integral, 
\begin{align*}
\uPi_x \ell_{\hat a(x)} & = \int_0^\infty \uPi_x(\{ \theta: \ell_{\hat a(x)}(\theta) > t \}) \, dt \\
& = \int_0^\infty \uPi_x(\{ \theta: \pi_x(\theta) < \beta_x(t) \}) \, dt \\
& = \int_0^\infty \beta_x(t) \, dt \\
& = \int_0^\infty \prior_x^\star(\{ \theta: \pi_x(\theta) < \beta_x(t) \}) \, dt \\
& = \int_0^\infty \prior_x(\{ \theta: \ell_{\hat a(x)}(\theta) > t \}) \, dt, \\
& = \prior_x^\star \ell_{\hat a(x)},
\end{align*}
where the second equality is by \eqref{eq:increasing2}, the third is by definition of $\uPi_x$ as a supremum of $\pi_x$, the fourth is by \eqref{eq:Qstar.equal}, the fifth is by \eqref{eq:increasing2}, and the sixth by the familiar formula for expectations of non-negative random variables.  This proves the second equality in \eqref{eq:risk}.  For the first equality, recall that $\uPi_x \ell_a \geq \prior_x^\star \ell_a$ for all $a$.  So, if there was another $a$ such that $\uPi_x \ell_a < \uPi_x \ell_{\hat a(x)}$, then that contradicts the definition \eqref{eq:bayes.rule} of $\hat a(x)$ as the minimizer. 
\end{proof}

This result generalizes the conclusion drawn in Example~\ref{ex:location} above; see also Example~\ref{ex:mises} below.  It can also be compared with the main result in \citet{taraldsen.lindqvist.2013} establishing that fiducial methods, which inherently have certain desirable confidence properties, also tend to produce good decision procedures.  Theorem~\ref{thm:risk} goes even further in the sense that IMs offer stronger and more comprehensive exact frequentist validity properties---i.e., they're not afflicted by false confidence the way fiducial solutions are---and they yield optimal decision procedures in some cases.  The additional $\preceq_x$-monotonicity condition \eqref{eq:increasing} on the loss in Theorem~\ref{thm:risk} complicates this comparison, but this isn't necessary for finding the IM's suggested action.  There's no reason to doubt the quality of the IM's optimal action when monotonicity fails, but this deserves further investigation. Indeed, our conjecture is that, for an invariant loss function that's not $\preceq_x$-monotone in the sense of \eqref{eq:increasing}, the minimizer of $a \mapsto \uPi_x\ell_a$ is still the Bayes rule $\hat a(x)$, but that $\uPi_x\ell_{\hat a(x)} > \prior_x^\star \ell_{\hat a(x)}$; see Example~\ref{ex:weibull} below.  If so, then, for example, when the right Haar-prior Bayes rule is admissible, minimax, minimum risk equivariant, etc., the possibilistic IM's upper risk minimizer is too.


\subsection{Examples}

\begin{figure}[t]
\begin{center}
\scalebox{0.6}{\includegraphics{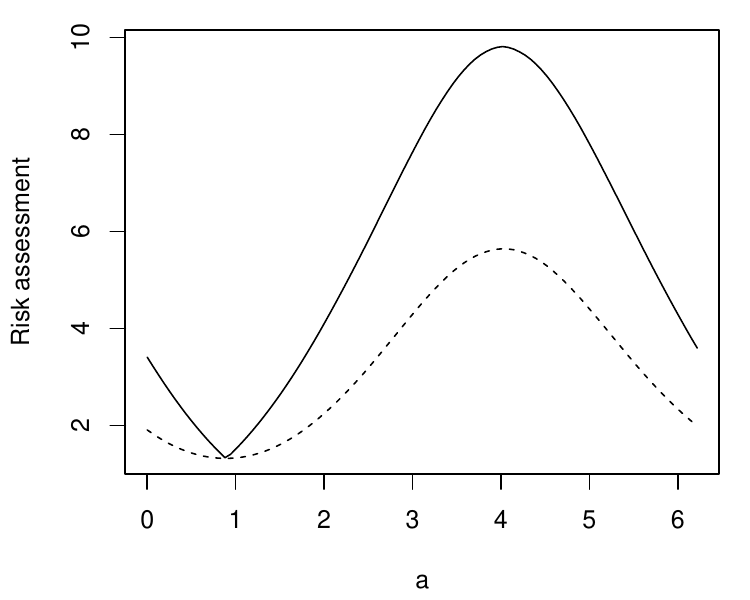}}
\end{center}
\caption{Plots of the Bayes risk (dashed) and the IM risk (solid) for the directional data illustration in Example~1.1 of \citet{mardia.jupp.book}, where the action $a$ is expressed in terms of an angle in the interval $[0,2\pi)$.}
\label{fig:vm.risk}
\end{figure}

\begin{ex}
\label{ex:mises}
As a quick illustration, consider a directional statistics application, as discussed in \citet{mardia.jupp.book}, where the data are supported on the unit circle.  A common model for such data is the so-called {\em von Mises distribution} with density function $p_\theta(x) \propto \exp\{ \kappa \, \cos(x-\theta) \}$, where $x$ and $\theta$ take values on the unit circle, and $\kappa > 0$ is a concentration parameter assumed here to be known.  This model is equivariant with respect to the group $\G$ of rotations, and the corresponding right Haar measure is the uniform distribution on the circle.  \citet{martin.isipta2023} used this as an illustration of Theorem~\ref{thm:fiducial} above, showing empirically---using the data in \citet[][Ex.~1.1]{mardia.jupp.book}---that the maximal element in the IM's credal set is indeed the Bayes posterior distribution based on the right Haar prior.  Here we follow up on that illustration to verify the new result in Theorem~\ref{thm:risk}.  Using that same data, Figure~\ref{fig:vm.risk} plots the Bayes risk and the IM's upper risk, as a function of the action $a$, based on the invariant loss function $\ell_a(\theta) = \cos^{-1}(\langle a, \theta \rangle)$, for $a$ in the unit circle; the plot transforms the unit circle to angles in the interval $[0,2\pi)$.  As Theorem~\ref{thm:risk} predicts, the two curves are not the same, but they agree in terms of both the location and value of their common minimum.  
\end{ex}

\begin{figure}[t]
\begin{center}
\subfigure[Model fit]{\scalebox{0.6}{\includegraphics{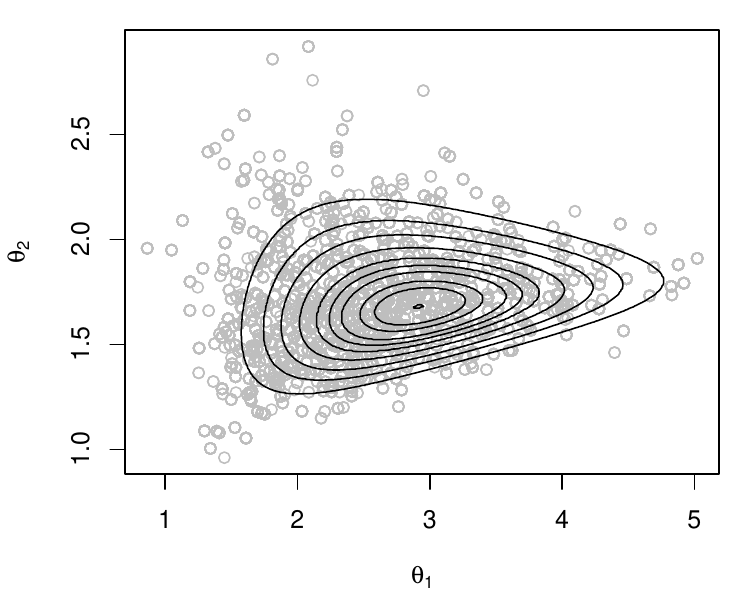}}}
\subfigure[Risk assessment]{\scalebox{0.6}{\includegraphics{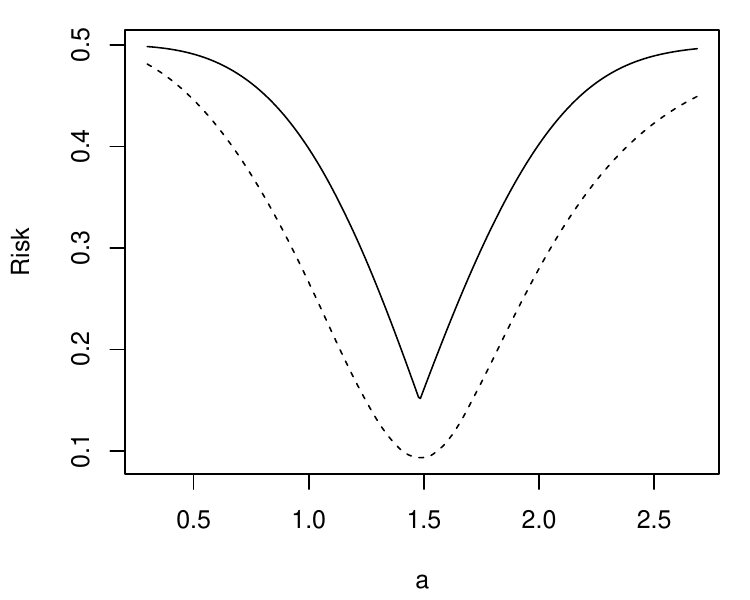}}}
\end{center}
\caption{Results for the Weibull data analysis in Example~\ref{ex:weibull}.  Panel~(a) plots the posterior samples (gray) and the possibilistic IM contour (black); Panel~(b) shows the Bayes risk (dashed) and the IM risk (solid).}
\label{fig:weib.risk}
\end{figure}

\begin{ex}
\label{ex:weibull}
Suppose that $X=(X_1,\ldots,X_n)$ consists of iid samples from a Weibull distribution with density function 
\[ p_\theta(x) = \frac{\theta_1 x^{\theta_1 - 1}}{\theta_2^{\theta_1}} \exp\Bigl\{ -\Bigl( \frac{x}{\theta_2} \Bigr)^{\theta_1} \Bigr\}, \quad x > 0. \]
Here $\theta_1$ is the shape parameter and $\theta_2$ the scale parameter, both are positive.  The Weibull model is equivariant with respect to the group of transformations $\G=\{g_{a,b}\}$ given by 
\[ g_{a,b}x := b x^{1/a}, \quad a > 0, \; b > 0. \]
That is, if $X_1 \sim {\sf Weib}(\theta_1, \theta_2)$, then $g_{a,b} X_1 \sim {\sf Weib}(a\theta_1, b\theta_2)$; hence $\bar g_{a,b}\theta = (a\theta_1, b\theta_2)$.  This transformation turns out to be just a location--scale transformation on the $\log X_1$ scale, with $\log\theta_2$ the location parameter and $\theta_2$ as the scale parameter.  So, one can take the right Haar prior for the location--scale formulation and transform it back to the shape--scale formulation, which gives $q(\theta) \propto (\theta_1\theta_2)^{-1}$; this agrees with the reference prior of \citet{sun.berger.1994} as summarized in \citet[][Sec.~39]{yang.berger.1998}.  For this illustration, we consider data of size $n=12$ presented in Problem~7.11 in \citet{ghosh-etal-book}, where it's suggested that the data be modeled by a Weibull distribution.  A relevant question concerns prediction of a future Weibull observation, and an informal/naive approach might be to use the mean of the data, $\bar X=1.49$.  More formally, if we treat the prediction as a generic action $a > 0$, then an invariant loss function is 
\[ \ell_a(\theta) = \bigl| P_\theta(a) - \tfrac12 \bigr|, \]
where $P_\theta$ is the Weibull cumulative distribution function.  Figure~\ref{fig:weib.risk} shows the Bayes and IM risk functions and, while the two curves are similar, they're clearly not the same.  But consistent with the conjecture stated at the very end of Section~\ref{SS:invariant.risk} above, the two risk functions have the same minimizer, $\hat a(x)=1.48$, which is the Bayes rule.  
\end{ex}

\section{Conclusion}
\label{S:discuss}

To date, IMs have only been used for inference; formal decision-making has not yet been investigated.  This paper fills this gap with roughly three contributions.  First, to motivate our proposed approach, we evaluated the behavior of the familiar Bayes assessment with respect to that of an oracle who knows the true $\Theta$ and determines the realized cost to the statistician for taking a particular action $a$.  We find that the Bayes assessment is overconfident in the sense that it tends to be smaller than the oracle assessment for some actions $a$, which puts the Bayesian at risk of systematically losing money even when he's taking gambles he himself deems to be acceptable.  Motivated by this shortcoming of probabilistic risk assessments, our second contribution is developing a framework for decision-making under uncertainty based on a possibilistic IM.  We define the IM's risk assessment, $a \mapsto \uPi_x\ell_a$, to be an upper expectation or Choquet integral of the loss with respect to the IM's possibility measure output $\uPi_x$, so that different actions $a$ can be compared using a natural variation of the minimize-expected-loss principle.  To address the Bayesian assessment's aforementioned overconfidence issue, we show that $\uPi_X \ell_a$ tends to not be smaller than the oracle assessment, uniformly over actions $a$, thus giving the IM a sense of reliability in the decision-making context which, in turn, provides some protection to the decision-maker.  Plus, a large-sample result establishes that the IM's decision-making safety isn't a consequence of making grossly conservative assessments.  Finally, the developments in Section~\ref{S:group} reveal that minimizing the IM's risk assessment yields actions that are high-quality, sometimes optimal.  Lots of numerical results are presented to back up the conceptual and theoretical developments mentioned above.  



There are a number of directions for future research that can be considered.  One is to 
investigate, both theoretically and empirically, the operating characteristics, such as frequentist risk, of the proposed IM-based decision rules in specific models.  Another is to determine if and how the monotonicity constraint in Theorem~\ref{thm:risk} can be relaxed.

An important notion in frequentist statistical decision theory is {\em admissibility}.  
Unfortunately, the IM's recommended action, defined in \eqref{eq:best}, is not admissible in general.  To see this, consider the $d$-dimensional normal mean problem.  Under loss $\ell_a(\theta) = \|a - \theta\|^2$, the IM's recommended action is $\hat a(x) = x$, which corresponds to the least squares, maximum likelihood, fiducial, and flat-prior Bayes estimator, was shown in, e.g., \citep{stein1956, brown1971} to be inadmissible when $d \geq 3$.  This is not especially surprising since the vacuous-prior version of the validity condition requires that the IM be ``unbiased'' in a certain sense which, in the present context, implies that the IM's recommended action be an unbiased estimator in the usual sense.  But the familiar bias--variance tradeoff suggests that the frequentist risk can be reduced by choosing a biased $\hat a$, e.g., like a proper-prior Bayes rule that's biased and admissible.  To introduce the appropriate ``bias'' into the IM construction, it's natural to consider incorporation of suitable (partial) prior information, e.g., sparsity.  Of course, incorporating such regularization would surely ruin the validity of the present IM construction, so one needs a modified version of validity, one that accounts for the knowledge encoded by the partial prior.  Such a notion of validity under partial prior information is being developed \citep{martin.partial2}, and an extension of the decision-theoretic results presented here to that case, along with admissibility considerations, are part of our ongoing work.

\section*{Acknowledgments}
The authors thank the anonymous reviewers for some very helpful suggestions on a previous version of this manuscript.  This work is partially supported by the U.S.~National Science Foundation, under grant DMS--2412628.

\appendix

\section{Further technical details}

\subsection{Computation}


While the IM construction is conceptually simple and its properties are strong, computation can be a challenge.  The go-to strategy has been to evaluate the contour function pointwise using the following simple/naive approximation, 
\[ \pi_{x}(\theta) \approx \frac1M \sum_{m=1}^M 1\{ R(X_{m,\theta}, \theta) \leq R(x, \theta) \}, \quad \theta \in \TT, \]
where $X_{m,\theta}$ are independent copies of the data $X$, drawn from $\prob_\theta$, for $m=1,\ldots,M$.  This is feasible at a few different $\theta$ values, but, e.g., evaluation over a fine grid covering the relevant portion of the parameter space $\TT$ can be expensive.  Recent developments \citep{immc} have revealed new opportunities for much more efficient computation. 

The strategies devised in \citet{immc} is based on the characterization of a credal set that corresponds to each $\alpha$ level of a possibilistic IM contour, as defined in \eqref{eq:credal}. Roughly, we can find the ``best'' probability distribution from the credal set to approximate the possibilistic IM. Since the probability measures in the credal set dominates the IM contour, the goal is to find the \emph{inner probabilistic approximation}---or, equivalently, the \emph{maximally diffuse element} in the possibilistic
IM’s credal set, as described in Section~\ref{SS:connection}---that exactly matches the $\alpha$ level, such that 
\[ 
\prior_x^{\text{inn}}\{ C_\alpha(x) \} = 1 - \alpha \quad \text{for each $\alpha \in [0, 1]$}.
\]
To further increase the computational efficiency, a second-layer approximation finds an \emph{outer approximation} of the inner approximation, which can be more conservative but much faster to compute. The implementation of this second-layer approximation assumes an elliptical form of each $\alpha$-cut and applies a Gaussian approximation that is stretched or shrunken to just barely encompass the inner approximation. The amount of stretching and shrinking varies for each $\alpha$. Samples of $\theta$ at each $\alpha$ level are obtained from this second-layer approximation and stitched together. In short, each $\alpha$-cut is assumed to be approximately ellipsoidal---justified by the asymptotic normality result in \citet{imbvm.ext}---but the entire sample is not approximately Gaussian but more like a mixture of Gaussian. Finally, a probability-to-possibility transform is introduced to use these probabilistic approximation samples to find a stitched possibilistic IM contour. 

Though the Monte Carlo approximation strategies, summarized above and detailed in \citet{immc}, have focused on finding the possibilistic IM contour, mainly useful for inference, they can be easily extended to risk assessment. One could naively use the probabilistic approximated samples to directly calculate the expected loss, but we find this strategy to under-estimate the upper expected loss. Alternatively, the stitched contour from the MC approximation can be plugged into the Choquet integral calculation, as defined in \eqref{eq:choquet}, saving a great amount of computation. 

\begin{figure}[t]
    \centering
    \subfigure[Comparison of contours]{\includegraphics[width=.49\linewidth]{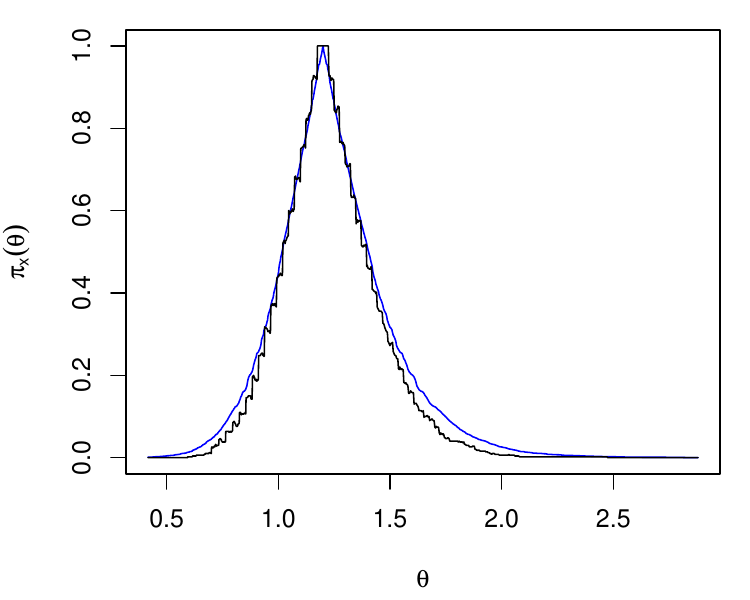}}
    \subfigure[Comparison of losses]{\includegraphics[width=.49\linewidth]{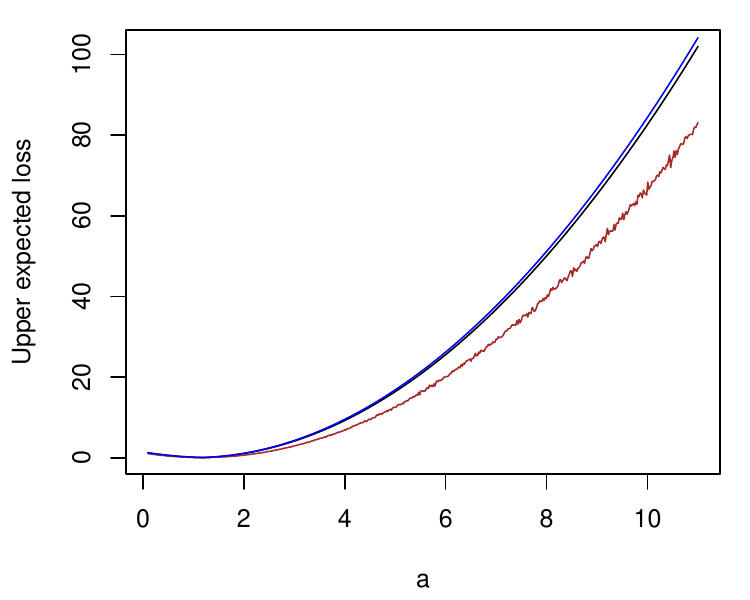}}
    \caption{Demonstration of Monte Carlo apprpoximation strategies. Panel (a) Compares exact (black) and stitched (blue) possibilistic IM contours. Panel (b) compares the upper expected loss from exact contour (black), stitched contour (blue), and Monte Carlo samples (brown).}
    \label{fig:MC_choquet}
\end{figure}

Figure~\ref{fig:MC_choquet} demonstrates empirical results from an example of a Poisson distribution with $n=20$. Panel (a) shows the contour function derived from a probability-to-possibility transform using Monte Carlo samples (blue) is quite close to the exact contour (black). Panel (b) shows that the loss directly calculated from approximated samples (brown) can be too low, but the approximate Choquet integral from the stitched IM contour (blue) is much closer to that calculated with the exact contour (black). In short, as long as the assumption of elliptical form is met (i.e. the model is assumed approximately Gaussian or the sample size is sufficiently large), the approximate Choquet integral should be sufficiently accurate.

\subsection{Large-sample results from Section~\ref{SS:asymptotics}} 
\label{AA:limit}

\subsubsection{Conditions}

\newcommand{\sL}{\mathcal{L}}

Certain regularity conditions are required in order to establish asymptotic concentration properties of estimators, posterior distributions, etc., and the same is true for IMs.  Here we briefly describe the so-called {\em Le Cam conditions} for asymptotic Gaussianity.  

We have a model $\{\prob_\theta: \theta \in \TT\}$ of probability distributions, supported on $\XX$, indexed by $\TT \subseteq \RR^D$, with $\prob_\theta$ having a density $p_\theta(x)$ relative to a $\sigma$-finite measure $\nu$ on $\XX$.  Then the data $X^n=(X_1,\ldots,X_n)$ is assumed to be independent and identically distributed of size $n$ with common distribution $\prob_\Theta$.  Following \citet[][Ch.~2]{bickel1998}, define the square-root density function as
\[ q_\theta(x) = p_\theta(x)^{1/2}. \]
The ``dot'' notation, e.g., $\dot g_\theta(x)$, represents a function that behaves like the derivative of $g_\theta(x)$ with respect to $\theta$ for pointwise in $x$.  If the usual partial derivative of $g_\theta(x)$ with respect to $\theta$ exists, then $\dot g_\theta(x)$ is that derivative; but suitable functions $\dot g_\theta(x)$ may exist even when the ordinary derivative fails to exist.  Finally, let $\sL_2(\nu)$ denote the set of measurable functions on $\XX$ that are square $\nu$-integrable. 

\begin{asmp}
The parameter space $\TT$ is open and there exists a vector $\dot q_\theta(x) = \{ \dot q_{\theta,d}(x): d=1,\ldots,D\}$, whose coordinates $\dot g_{\theta,d}$ are elements of $\sL_2(\nu)$, such that the following conditions hold: 
\begin{itemize}
\item[R1.] the maps $\theta \mapsto \dot q_{\theta,d}$ from $\TT$ to $\sL_2(\nu)$ are continuous for each $d=1,\ldots,D$; 
\item[R2.] at each $\theta \in \TT$, 
\begin{equation}
\label{eq:dqm}
\int \bigl| q_{\theta + u}(x) - q_\theta(x) - u^\top \dot q_\theta(x) \bigr|^2 \, \nu(dx) = o(\|u\|^2), \quad u \to 0 \in \RR^D; 
\end{equation}
\item[R3.] and the $D \times D$ matrix $\int \dot q_\theta(x) \, \dot q_\theta(x)^\top \, \nu(dx)$ is non-singular for each $\theta \in \TT$. 
\end{itemize} 
\end{asmp}

The condition in \eqref{eq:dqm} is often described as $\theta \mapsto q_\theta$ being {\em differentiable in quadratic mean}.  Note that this condition does not require the square-root density to actually be differentiable at $\theta$, only that it be ``locally linear'' in a certain average sense.  The classical Cram\'er conditions assume more than two continuous derivatives, so \eqref{eq:dqm}, which does not even require existence of a first derivative, is significantly weaker; sufficient conditions for \eqref{eq:dqm} are given in \citet[][Lemma~7.6]{vaart1998}.  Then the {\em score function} $s_\theta$ is defined in terms of $\dot q_\theta$ as 
\[ s_\theta(x) = \frac{2\dot q_\theta(x)}{q_\theta(x)} \, 1\{q_\theta(x) > 0\}, \]
and it can be shown that $\int s_\theta(x) \, \prob_\theta(dx) = 0$ for each $\theta$.  Moreover, Condition~R3 above implies non-singularity of the Fisher information matrix $I_\theta = \int s_\theta(x) s_\theta(x)^\top \, \prob_\theta(dx)$ for each $\theta \in \TT$.  \citet[][Prop.~2.1.1]{bickel1998} provide sufficient conditions for R1--R3.   

One further condition is required by \citet{imbvm.ext} is that the maximum likelihood estimator is consistent, i.e., $\hat\theta_{X^n} \to \Theta$ in $\prob_\Theta$-probability as $n \to \infty$.  This, of course, is not automatic, but holds quite broadly when the dimension of the parameter space is fixed, as we are assuming here.

\subsubsection{Notation and background}

To set the scene for proof the follows, for a mean vector $\mu \in \RR^D$ and a covariance matrix $\Sigma \in \RR_+^{D \times D}$, define the Gaussian possibility contour $\gamma_{\mu, \Sigma}$ as 
\begin{align*}
\gamma_{\mu, \Sigma}(z) & := \prob\bigl[ \exp\{ -\tfrac12 (Z-\mu)^\top \Sigma^{-1} (Z-\mu) \} \leq \exp\{ -\tfrac12 (z-\mu)^\top \Sigma^{-1} (z-\mu) \} \bigr] \\
& = \prob\bigl\{ (Z-\mu)^\top \Sigma^{-1} (Z-\mu) \geq (z-\mu)^\top \Sigma^{-1} (z-\mu) \bigr\} \\
& = 1 - F_D\bigl\{ (z-\mu)^\top \Sigma^{-1} (z-\mu) \bigr\}, 
\end{align*}
where $F_D$ is the distribution function of the $\chisq(D)$ distribution.  Other authors have used a different form of Gaussian possibility \citep[e.g.,][]{denoeux.fuzzy.2022, denoeux.fuzzy.2023}, but we find the one above appealing and convenient because the Gaussian probability distribution $\nm_D(\mu, \Sigma)$ is the maximally diffuse element of the credal set determined by the contour $\gamma_{\mu, \Sigma}$.   

Of particular relevance to our analysis here is the case 
\begin{equation}
\label{eq:gauss.countour}
\gamma_{X^n}(\theta) := \gamma_{\Theta + n^{-1/2} \Delta_\Theta(X^n), (n I_\Theta)^{-1}}(\theta),
\end{equation}
i.e., the Gaussian possibility contour with (data- and true-parameter-dependent) mean vector $\mu = \Theta + n^{-1/2} \Delta_\Theta(X^n)$, where
\[ \Delta_\theta(X^n) := n^{-1/2} I_\theta^{-1} \sum_{i=1}^{n}  \frac{\partial}{\partial\theta} \log p_\theta(X_i), \]
and covariance matrix $\Sigma = (nI_\Theta)^{-1}$, where, recall $I_\theta$ is the Fisher information matrix derived from the model.  A relevant detail both in \citet{imbvm.ext} and in what follows is that $\Delta_\Theta(X^n) \to \nm_D(0, I_\Theta^{-1})$ in distribution, as $n \to \infty$, under $\prob_\Theta$.  Of course, we don't/can't use the contour $\gamma_{X^n}$ in practice, since we don't know $\Theta$; this is only for the theoretical analysis.  Then \citet{imbvm.ext} showed that, for any compact set $\K \subset \TT$, 
\[ \sup_{\theta \in \K} \bigl| \pi_{X^n}(\theta) - \gamma_{X^n}(\theta) \bigr| \to 0, \quad \text{in $\prob_\Theta$-probability as $n \to \infty$}. \]
From here, a number of different properties can be established concerning the upper probabilities assigned to set $H \subseteq \TT$ and to upper expectations/Choquet integrals.  The latter case is relevant in what comes next.

\subsubsection{Proof of Theorem~\ref{thm:limit}}

\newcommand{\ellsupp}{\mathscr{L}}

Write $\uGamma_{X^n}$ for the possibility measure determined by the contour function $\gamma_{X^n}$ in \eqref{eq:gauss.countour}.  The strategy of the proof is to first demonstrate that $\uPi_{X^n}$ and $\uGamma_{X^n}$ merge asymptotically relative to their upper expected losses, and then to show that the upper expected loss under $\uGamma_{X^n}$ merges with the oracle.  

As before, write $\uGamma_{X^n} \ell_a$ for the upper expected loss, defined via a Choquet integral with respect to $\uGamma_{X^n}$; see below.  Then a simple application of the triangle inequality gives
\begin{equation}
\label{eq:main.triangle}
|\uPi_x \ell_a - \ell_a(\Theta)| \le |\uPi_x \ell_a - \uGamma_x \ell_a| + |\uGamma_x \ell_a - \ell_a(\Theta)|.
\end{equation}
We will show separately that both terms on the right side vanish asymptotically.

There are different but equivalent expressions for the Choquet integral with respect to a possibility measure; see, e.g., Proposition~7.14 in \citet{lower.previsions.book}.  We've been using one of those versions in this paper up to now, but a different version---with the positions of the integrand and the contour swapped---is more convenient for the analysis that follows.  In particular, we work the upper expected losses 
\begin{align*}
\uPi_{X^n}\ell_a & = \inf \ell_a + \int_{\inf \ell_a}^{\sup \ell_a} \sup_{\theta: \ell_a(\theta) \geq t} \pi_{X^n}(\theta) \, dt \\
\uGamma_{X^n}\ell_a & = \inf \ell_a + \int_{\inf \ell_a}^{\sup \ell_a} \sup_{\theta: \ell_a(\theta) \geq t} \gamma_{X^n}(\theta) \, dt.
\end{align*}
To simplify notation, we'll assume without loss of generality that the minimum loss is 0, so that the additive term ``$\inf \ell_a$'' above can be dropped.  

For the first term in \eqref{eq:main.triangle}, 
\begin{align*}
|\uPi_{X^n} \ell_a - \uGamma_{X^n} \ell_a| & = \Bigl| \int_{0}^{\sup\ell_a}\sup_{\theta : \ell_a(\theta) \ge t} \pi_{X^n}(\theta) \, dt - \int_{0}^{\sup\ell_a}\sup_{\theta : \ell_a(\theta) \ge t} \gamma_{X^n}(\theta) \, dt \Bigr| \\
& \le \Bigl| \int_{0}^{\sup\ell_a}\sup_{\theta : \ell_a(\theta) \ge t} \pi_{X^n}(\theta) - \sup_{\theta : \ell_a(\theta) \ge t} \gamma_{X^n}(\theta) \, dt \Bigr| \\
& \le \int_{0}^{\sup\ell_a}\sup_{\theta : \ell_a(\theta) \ge t} \bigl| \pi_{X^n}(\theta) - \gamma_{X^n}(\theta) \bigr| \, dt \\
& \le (\sup\ell_a) \times \sup_{\theta \in \ellsupp}|\pi_{X^n}(\theta) - \gamma_{X^n}(\theta)|,
\end{align*}
where $\ellsupp$ is the compact support of $\ell_{a}$.  Because $\ell_{a}$ is continuous with compact support, it follows that $\sup \ell_a < \infty$.  Then Theorem~1 in \citet{imbvm.ext} implies that the upper bound vanishes; hence, $|\uPi_x \ell_a - \uGamma_x \ell_a| \to 0$ in $\prob_{\Theta}$-probability as $n \to \infty$.

It remains to show that the second term on the right side of \eqref{eq:main.triangle} vanishes asymptotically.  First, since $\cred(\uGamma_x)$ contains a point mass distribution at $\Theta + n^{-1/2}\Delta_\Theta(X^n)$, it follows immediately that $\uGamma_x \ell_a \geq \ell_a(\Theta + n^{-1/2}\Delta_\Theta(X^n))$.  Since $\Delta_\Theta(X^n)$ is bounded in $\prob_\Theta$-probability, it follows from the continuous mapping theorem that $\ell_a(\Theta + n^{-1/2}\Delta_\Theta(X^n)) \to \ell_a(\Theta)$ in $\prob_\Theta$-probability.  Next, for any $\eps > 0$
\[ \uGamma_{X^n} \ell_a = \int \sup_{\theta: \ell_a(\theta) \geq t} \gamma_{X^n}(\theta) \, dt = \Bigl( \int_0^{\ell_a(\Theta) + \eps} + \int_{\ell_a(\Theta) + \eps}^{\sup \ell_a} \Bigr) \sup_{\theta: \ell_a(\theta) \geq t} \gamma_{X^n}(\theta) \, dt. \]
For the first integral, use the trivial bound $\sup \gamma_{X^n}(\theta) \leq 1$, and, for the second integral, use the fact that the integrand is non-increasing in $t$.  This gives the upper bound 
\[ \uGamma_{X^n} \ell_a \leq \{\ell_a(\Theta) + \eps\} + \{\sup \ell_a - \ell_a(\Theta)-\eps\} \sup_{\theta: \ell_a(\theta) \geq \ell_a(\Theta) + \eps} \gamma_{X^n}(\theta). \]
The set $\{\theta: \ell_a(\theta) > \ell_a(\Theta) + \eps\}$ is contained in the complement of an open neighborhood around $\Theta$, say, $\{\theta : \|\Theta - \theta\| < \delta_{\eps}\}$, where $\delta_\eps$ depends on $\eps$ and possibly on $a$.  Furthermore, since $\Delta_\Theta(X^n) \to \nm_{D}(0,I_{\Theta}^{-1})$ in distribution, and hence $n^{-1/2} \Delta_\Theta(X^n) \to 0$ in $\prob_\Theta$-probability, it follows that for sufficiently large $n$, with $\prob_{\Theta}$-probability $\to 1$,
\[
\{\theta  :  \|\Theta - \theta\| < \delta_{\eps}\} \supset B_{n} := \{ \theta  :  \|\Theta + n^{-1/2} \Delta_\Theta(X^n) - \theta\| < n^{-1/4}\},
\]
and so $\{\theta: \ell_a(\theta) \geq \ell_a(\Theta) + \eps\} \subset B_{n}^{c}$.  For any $\theta \in B_{n}^{c}$, the Gaussian contour in \eqref{eq:gauss.countour} is of the form $\gamma_{X^n}(\theta) = 1 - F_D(\square)$, where the argument $\square$ is
\begin{align*}
\square & = \{\Theta + n^{-1/2} \Delta_\Theta(X^n) - \theta\}^\top (nI_{\Theta})\{\Theta + n^{-1/2} \Delta_\Theta(X^n) - \theta\} \\
& \ge n \lambda_\text{min}(I_{\Theta}) \cdot \|\Theta + n^{-1/2} \Delta_\Theta(X^n) - \theta\|^{2} \\
& \ge n^{1/2} \lambda_\text{min}(I_{\Theta}),
\end{align*}
where $\lambda_\text{min}(A)$ denotes the minimum eigenvalue of $A$ and the last line follows because $\|\Theta + n^{-1/2} \Delta_\Theta(X^n) - \theta\|^{2} \geq n^{-1/4}$ for all $\theta \in B_n^c$.  Accordingly, $\gamma_{X^n}(\theta) \to 0$ in $\prob_\Theta$-probability uniformly for $\theta \in B_n^c$ and, therefore, 
\[ 
\sup_{\theta: \ell_a(\theta) > \ell_a(\Theta) + \eps} \gamma_{x}(\theta) \le \sup_{B_{n}^{c}} \gamma_{x}(\theta) \to 0 \quad \text{in $\prob_\Theta$-probability}. 
\]
This implies $\limsup_{n \to \infty} \uGamma_{X^n}\ell_a \leq \ell_a(\Theta) + \eps$.  Putting this together with the simple lower bound mentioned previously gives 
\[ \ell_a(\Theta) \leq \liminf_{n \to \infty} \uGamma_{X^n} \ell_a \leq \limsup_{n \to \infty} \uGamma_{X^n}\ell_a \leq \ell_a(\Theta) + \eps, \quad \text{in $\prob_\Theta$-probability}. \]
But since $\eps > 0$ is arbitrary, it must be that $\uGamma_{X^n}\ell_a \to \ell_a(\Theta)$ in $\prob_\Theta$-probability.  Combining this with the first part of the proof gives 
\[ \bigl| \uPi_{X^n} \ell_a - \ell_a(\Theta) \bigr| \to 0 \quad \text{in $\prob_\Theta$-probability as $n \to \infty$}. \]

Finally, for the ``uniform-in-actions'' version of the result, we only need to make two observations.  First, since $\{\ell_a: a \in \action\}$ is assumed to be bounded, the term ``$\sup \ell_a$'' in the bound of $|\uPi_{X^n} \ell_a - \uGamma_{X^n} \ell_a|$ actually doesn't depend on the action $a$.  Consequently, that term vanishes uniformly over $a \in \action$ in $\prob_\Theta$-probability.  Second, since $\{\ell_a: a \in \action\}$ is equicontinuous, the ``$\delta_\eps$'' in the bound of $|\uGamma_{X^n}\ell_a - \ell_a(\Theta)|$ also doesn't depend on $a$ and, likewise, the corresponding upper bound vanishes uniformly over $a \in \action$ in $\prob_\Theta$-probability.  Putting these pieces together gives the ``merges uniformly'' claim. 

\bibliographystyle{apalike}
\bibliography{mybib}

\end{document}